\documentclass[12pt]{article}
\usepackage{}
\usepackage{amssymb}
\usepackage{amsthm}
\usepackage[utf8]{inputenc}
\usepackage{extarrows}
\usepackage{cases}
\usepackage{graphicx}
\usepackage{mathrsfs}
\usepackage{graphics}
\usepackage[misc]{ifsym}
\usepackage{color}
\usepackage{indentfirst}
\usepackage{times}

\numberwithin{figure}{section}
\numberwithin{equation}{section}

\newtheorem{theorem}{Theorem}[section]
\newtheorem{proposition}[theorem]{Proposition}
\newtheorem{definition}[theorem]{Definition}
\newtheorem{corollary}[theorem]{Corollary}
\newtheorem{lemma}[theorem]{Lemma}
\newtheorem{remark}[theorem]{Remark}
\newtheorem{example}[theorem]{Example}

\newcommand{\cA}{{\mathcal A}}
\newcommand{\cB}{{\mathcal B}}

\newcommand{\cD}{{\mathcal D}}

\newcommand{\cM}{{\mathcal M}}
\newcommand{\cN}{{\mathcal N}}

\newcommand{\sM}{{\mathscr M}}

\def\T{\mathbb{T}}
\def\Z{\mathbb{Z}}
\def\R{\mathbb{R}}

\def\N{\mathbb{N}}

\def\1{\mathbb{1}}

\def\bc{\begin{center}}
\def\ec{\end{center}}
\def\be{\begin{equation}}
\def\ee{\end{equation}}
\def\ba{\begin{array}}
\def\ea{\end{array}}
\def\benu{\begin{enumerate}}
\def\eenu{\end{enumerate}}
\def\bt{\begin{theorem}}
\def\et{\end{theorem}}
\def\bl{\begin{lemma}}
\def\el{\end{lemma}}
\def\bco{\begin{corollary}}
\def\eco{\end{corollary}}
\def\bn{\begin{numcases}}
\def\en{\end{numcases}}
\def\br{\begin{remark}}
\def\er{\end{remark}}
\def\bd{\begin{definition}}
\def\ed{\end{definition}}
\def\bp{\begin{proposition}}
\def\ep{\end{proposition}}
\def\bo{\begin{proof}}
\def\eo{\end{proof}}
\def\bx{\begin{example}}
\def\ex{\end{example}}

\def\a{\alpha}
 \def\de{\delta}

\def\lam{\lambda} 

\def\ve{\varepsilon}
\def\sig{\sigma}\def\Sig{\Sigma}
\def\vp{\varphi}
\def\w{\omega}

\def\~{\widetilde}
\def\A{\forall}

\def\Cup{\bigcup}

\def\ra{\rightarrow}

\def\8{\infty}
\def\X{\times}

\def\mb{\mbox}

\def\ss{\subset}

\def\.{\cdot}\def\leq{\leqslant}\def\geq{\geqslant}\def\d{{\rm d}}

\def\Hs{\hspace{0.8cm}}
\def\hs{\hspace{0.4cm}}
\def\Vs{\vskip10pt}
\def\vs{\vskip5pt}

\def\[{\left[}
\def\]{\right]}
\def\({\left(}
\def\){\right)}

\title{
\Vs
 On the forward dynamical behavior of nonautonomous lattice dynamical systems
  \footnote{This work is supported by NSF of China under Grants
 11871368, 11801190.}}

\author{\small
{Chunqiu Li,}\footnote{Corresponding author. {\em E-mail}: licqmath@tju.edu.cn}
\Hs {Jintao Wang} \footnote{ Email: wangjt@wzu.edu.cn}
\\
\small\it Department of Mathematics, Wenzhou University,\\ \small\it Wenzhou, Zhejiang Province, 325035, P. R. China\\
}
\date{}
\begin{document}
\maketitle
\begin{abstract}\baselineskip 14pt
In this article, we study the forward dynamical behavior of nonautonomous lattice systems.
We first construct a family of sets $\{\cA_\ve(\sig)\}_{\sig\in \Sig}$
in arbitrary small neighborhood of a global attractor of the skew-product flow generated by  a general nonautonomous lattice system,
which is forward invariant and uniformly forward attracts any bounded subset of the phase space.
Moreover, under some suitable conditions, we further construct a family of sets $\{\cB_\ve(\sig)\}_{\sig\in \Sig}$ such that it uniformly  forward exponentially attracts bounded subsets of the phase space.
As an application, we study the discrete Gray-Scott model in detail and illustrate
how to apply our abstract results to some concrete lattice system.

\vs
 \vs
\end{abstract}


\setcounter {equation}{0}
\section{Introduction}
It is well known that the global attractor is an appropriate object to describe the asymptotic behavior of autonomous systems, and the corresponding theory of attractors has been fully developed in the past decades; see \cite{Hale,LWX17,SY,T}. In contrast, the case for nonautonomous systems seems to be far more difficult, though the study of the  dynamics of nonautonomous systems has been paid much attention in recent years; see e.g. \cite{CLR,Chep}. In fact, even  the notions of attractors are still under investigations.


Similar to autonomous systems, there are several concepts, uniform attractors, kernel sections, pullback attractors and forward attractors  which can be used to capture the dynamics of nonautonomous dynamical systems (NDSs in short).
These attractors can be regarded as a natural generalization of global attractors of autonomous systems. Among these attractors, pullback attractors and forward attractors are very important and extensively studied by researchers.
On the existence of a pullback attractor, there are quite general results for both lattice systems and continuous systems; see e.g. \cite{Cara,CLR,FW,WZZ,XL21,ZXL}, etc. On the other hand, in most cases we are  more concerned with the evolution of a nonautonomous lattice system in the future. However, some relative works mainly concentrate on the continuous case and its pullback asymptotic behavior; see e.g. \cite{CLO,JLD,JLLQ}. In general, a pullback attractor may not be forward attracting and provide  little  information on this aspect.

Lattice dynamical systems (LDSs in short) are spatiotemporal systems
with discretization in some variables, which appear in many
different fields such as biology \cite{K}, chemical reaction theory
\cite{EN}, electrical engineering \cite{CP} and so on. Nowadays, the asymptotic behavior of LDSs has attracted much attention from researchers; see e.g. \cite{A11,BLW,CMV,HK,KY,LSLZ,W1,Zhou1,Zhou2,ZZL}.
Particularly, Zhou and Zhao \cite{ZZL} gave some sufficient and necessary conditions for the existence of a compact uniform attractor  in the framework of processes on a Hilbert space of infinite sequences. Concerning the uniform attractors, kernel
sections, pullback attractors of LDSs,  one can find a vast body of literature on their existence and structure; see e.g. \cite{W,ZZ07,ZZ,ZZ2,ZZ3,ZH}.

In this article, we are concerned with the forward dynamical behavior of general nonautonomous LDSs. Generally speaking, the forward attractor may be quite  suitable  to capture the forward dynamical behavior of a nonautonomous LDS.
Unfortunately, there are few results on the existence  of such type of attractors except in some particular cases of continuous systems such as the periodic and the asymptotically autonomous ones, and these problems are still undergoing investigations; see e.g. Carvalho {\it et al.} \cite[p. 595]{CLRS}, Cui {\it et al.} \cite{CK}, Cheban {\it et al.} \cite{CKS}, Kloeden {\it et al.} \cite{KL}, Langa {\it et al.} \cite{LOS}, and Wang {\it et al.} \cite{WLK}. It is well known that if the pullback attraction of a pullback attractor $\cA=\{A(\sig)\}_{\sig\in \Sig}$ of a family of processes $\{U_\sig(t,0)\}_{t\geq 0},\sig\in \Sig$ in the space $E$ is uniform with respect to $\sig\in \Sig$, where $\Sig$ is a compact metric space and $E$ is a complete metric space, then it is also a forward attractor. Using this result,  the authors \cite{CKS} considered the continuous system and verified that  if the section  $A(\sig)$ of $\cA$ is lower semicontinuous in $\sig$, then $\cA$ is a forward attractor. However, it is hard to prove the lower semicontinuity of a family of sets generally.

Motivated by these works mentioned above, we shall construct a family of sets, which may be used to describe the forward dynamical behavior of a general LDS in some sense. Specifically, let $\{U_\sig(t,0)\}_{t\geq 0},\sig\in \Sig$ be a family of processes generated by the LDS in a phase space $E$, where $\sig\in\Sig$ and $\Sig$ is called a symbol space. If the mappings $(u,\sig)\ra U_\sig(t,0)u$ are continuous from $E\X\Sig$ to $E$, the family of continuous processes can reduce to a autonomous semigroup $\Psi$, called a skew-product flow. It will be  shown that in any small $\ve$-neighborhood of the global attractor $\cA$ of $\Psi$, there exists a family of sets $\{\cA_\ve(\sig)\}_{\sig\in \Sig}$, which is forward invariant for $\{U_\sig(t,0)\}_{t\geq 0},\sig\in \Sig$, such that it uniformly forward attracts each bounded subset of $E$. 

The second result of this article is to construct a family of sets $\{\cB_\ve(\sig)\}_{\sig\in \Sig}$ in an arbitrary small $\ve$-neighborhood of the global attractor such that it uniformly forward exponentially attracts each bounded subset of $E$. For this purpose, we first construct a pullback exponential attractor for $\{U_\sig(t,0)\}_{t\geq 0},\sig\in \Sig$. 
As we all know, an exponential attractor of an autonomous system, which was first constructed by Eden {\it et al.} \cite{EFNT}, attracts the trajectories at an exponential rate and is more robust under perturbations.  For nonautonomous systems, pullback exponential attractors are extensively studied; see e.g. \cite{ZH} for lattice systems and \cite{CS,CS1,EMZ,EYY,Zhong} for continuous systems.
Particularly, the authors \cite{EYY} constructed pullback exponential attractors of nonautonomous systems in the framework of a process on Banach spaces. Later, based on the results in \cite{EYY}, Zhou and Han \cite{ZH} presented a sufficient and necessary condition for the existence of pullback exponential attractors of the process generated by a nonautonomous LDS. The interested readers are referred to \cite{A,BN,DN,FY,Han} for some concrete continuous and lattice systems on the existence and construction of exponential attractors. By the exponential attractivity of the pullback exponential attractor, we will further show that the family of sets $\{\cB_\ve(\sig)\}_{\sig\in\Sig}$ uniformly forward exponentially attracts each bounded subset of the phase space under some suitable conditions.

As an application of our main results, we consider the nonautonomous discrete three-component reversible Gray-Scott model driven by quasi-periodic external forces.
When there are no external forces, the Gray-Scott equations reduce to
the corresponding autonomous Gray-Scott model, which was first
introduced by H. Mahara {\it et al.} \cite{M} to describe the isothermal,
cubic autocatalytic continuously fed and diffusive reactions of two
chemicals \cite{Y2}. Concerning
the initial boundary value and Cauchy problems for both discrete and continuous
Gray-Scott equations, there are many works in its
autonomous and nonautonomous cases; see \cite{JZY,Y2,Y3} and
the references therein.

The paper is organized as follows. In section 2, we make some preliminaries. Section 3 is devoted to our main results. In section 4 we give an example to illustrate our results. Finally, we give a remark on some extensions.
\setcounter {equation}{0}

\section{Preliminaries}
In this section, we collect some basic concepts on dynamical systems. One can also see \cite{CLR,Chep} etc., for details.
\vs
Let $E$ be a complete metric space and $\Sig$ a metric space with the metrics $d$ and $\rho$, respectively. Given two subsets $A,B$ of $E$. Define the Hausdorff-semidistance and Hausdorff distance of $A$ and $B$, respectively, as $${\rm d}_H(A,B)=\sup_{x\in A}d(x,B),\hs {\rm \de}_H(A,B)=\max\{{\rm d}_H(A,B),{\rm d}_H(B,A)\}.$$ The $\ve$-neighborhood of $A$ is defined by the set $$N_\ve(A)=\{x\in X: d(x,A)<\ve\}.$$

Let $E$ be a Banach space with norm $\|\.\|$.  A function $f(t)\in L_{loc}^2(\R,E)$ is called {\it translation bounded} in $L_{loc}^2(\R,E)$ if $$\|f\|_{L_b^2(\R,E)}=\sup_{t\in\R}\int^{t+1}_{t}\|f\|_E^2{\rm d}s<+\8.$$ Denote $L_b^2(\R,E)$ the sets of translation bounded functions in $L_{loc}^2(\R,E)$.

A function $f(t)\in L_{loc}^2(\R,E)$ is said to be {\it translation compact} in $L_{loc}^2(\R,E)$, if the set $\{f(\tau+\cdot):\tau\in\R\}$ is precompact in $L_{loc}^2(\R,E)$. By $L_{c}^2(\R,E)$ we denote the union of translation compact functions
in $L_{loc}^2(\R,E)$. Then $L_{c}^2(\R,E)\ss L_{loc}^2(\R,E)$ (see \cite{LWZ}).

\vs
Let $E$ be a Banach space. A family of  mappings $\{U(t,\tau)\}_{t\geqslant\tau}$ is called {\it a process} in $E$ if
\benu
\item[(i)] $U(t,s)U(s,\tau)=U(t,\tau), \hs t\geqslant s\geqslant\tau, \tau\in \R;$
\item[(ii)] $U(\tau,\tau)={\rm Id}_E,\hs \tau\in\R.$
\eenu
$\{U_\sig(t,\tau)\}_{t\geqslant\tau},\sig\in \Xi$ is said to be {\it a family of processes} in $E$, if for any $\sig\in \Xi$, $\{U_\sig(t,\tau)\}_{t\geqslant\tau}$ is a process in $E$, where $\Xi$ is called {\it a symbol space}.
\bd
A set $B_0\ss E$ is said to be uniformly (w.r.t. $\sig\in\Xi$) absorbing for the family of processes $\{U_\sig(t,\tau)\}_{t\geqslant\tau},\sig\in \Xi$, if for each $\tau\in \R$ and bounded set $B\ss E$, there exists $t_0=t_0(\tau,B)\geqslant \tau$ such that $$\Cup_{\sig\in \Xi}U_\sig(t,\tau)B\ss B_0,\hs t\geqslant t_0.$$
\ed
\bd
A closed set $\cA_\Xi$ is said to be a uniformly (w.r.t. $\sig\in \Xi$) attractor for the family of processes $\{U_\sig(t,\tau)\}_{t\geqslant\tau},\sig\in \Xi$, if
\benu
\item[{\rm (i)}] $\lim\limits_{t\ra+\8}\sup\limits_{\sig\in\Xi}{\rm d}_H(U_\sig(t,\tau)B,\cA_\Xi)=0,$ for each fixed $\tau\in\R$ and each bounded set $B\ss E$;
\item[{\rm (ii)}] $\cA_\Xi$ is the minimal set (for inclusion relation) among those satisfying {\rm (i)}.
\eenu
\ed
Now we always assume that the following assumptions hold true.

\noindent Assumption {\bf (I)}.
Let $\{\theta_s: s\geq 0\}$ be the natural translation semigroup (see \cite{Chep}) acting on $\Xi$ and satisfy
\benu
\item[{\rm (a)}] $\theta_s\Xi=\Xi$ for all $s\geq 0$;
\item[{\rm (b)}] $U_\sig(t+s,\tau+s)=U_{\theta_s\sig}(t,\tau)$.
\eenu

\bl\label{2.4}
Let $\tau_0\in \R$ be fixed. Then for any $\tau\in \R$ and $\sig\in \Sig$, there exists unique $\sig_1\in \Sig$ such that
\be\label{e2.3}
U_\sig(t,\tau)=U_{\sig_1}(t-\tau+\tau_0,\tau_0),\hs t\geqslant \tau.\ee
\el
\bo
 Let $\sig\in \Sig$. Since $$U_\sig(t,\tau)=U_{\sig}(t-\tau+\tau_0+\tau-\tau_0,\tau-\tau_0+\tau_0),$$ we deduce from Assumption ({\bf I}) that \eqref{e2.3} holds true with $\sig_1=\theta_{\tau-\tau_0}\sig$. Note that $\theta_t:\Sig\mapsto \Sig$ is a homeomorphism for each $t\in \R$.  The result of the lemma holds.
\eo
\br\label{r2.6}
If $\tau_0=0,$ then the above lemma shows $$U_\sig(t,\tau)=U_{\sig_1}(t-\tau,0),\hs t\geqslant \tau$$ for each $\tau\in \R$ and $\sig\in \Sig$.
\er

Finally, we recall the definition of a pullback exponential attractor for the family of processes $\{U_\sig(t,0)\}_{t\geq 0},\sig\in \Sig.$
\bd
A family of sets $\{\sM(\sig)\}_{\sig\in\Sig}$ is said to be a pullback exponential attractor $\{U_\sig(t,0)\}_{t\geq 0},\sig\in \Sig$, if it satisfies the following properties:
\benu
\item[(1)] $\sM(\sig)$ is a compact subset of $E$ for all $\sig\in \Sig$, and its fractal dimension is finite uniformly with respect to $\sig\in \Sig$, i.e. $\sup_{\sig\in \Sig}\dim_f\sM(\sig)<+\8.$
\item[(2)] It is forward invariant: $$U_\sig(t,0)\sM(\sig)\subset \sM(\theta_t\sig),\Hs t\geqslant 0,\hs \sig\in \Sig.$$
\item[(3)] There exists $\a>0$ such that for each $\sig\in \Sig$ and any bounded subset $B$ of $E$, there exists $T_{\sig,B}>0$ with
$$ {\rm d}_H(U_\sig(0,-t)B,\sM(\sig))\leqslant C_{B}{\rm e}^{-\a t},\hs  \forall t\geqslant T_{\sig,B},$$ where $C_{B}$ is a positive constant depending on $B$.
\eenu
\ed


\setcounter {equation}{0}
\section{Forward dynamical behavior of LDSs}

In this section we consider the forward dynamical behavior of general nonautonomous LDSs.

Set
$$\ell^2=\{u=(u_{m})_{m\in \mathbb{Z}}:\  u_{m}\in\mathbb{R},\hs
  \sum\limits_{m\in\mathbb{Z}}u_{m}^{2}<+\infty\},
$$
which is equipped  with the inner product and norm defined by
$$
  (u,v)=\sum\limits_{m\in\mathbb{Z}}u_{m}v_{m}, \hs \|u\|^2=(u,u),\Hs \A u=(u_m)_{m\in \mathbb{Z}},\, v=(u_m)_{m\in \mathbb{Z}}\in \ell^{2}.$$
For convenience, we use $E$ to denote $\ell^2$. Then $E$ is a Hilbert space.

Consider the following abstract nonautonomous LDSs in $E$:
\begin{align}
&\dot {u}=f(u)+\sig(t),\hs t>t_0;\label{e3.1}\\
&u(t_0)=u_0,\hs t_0\in\R, \label{e3.2}
\end{align}
where $u=(u_m)_{m\in\Z}, u_0=(u_{m,0})_{m\in\Z}\in E$, $f=(f_m)_{m\in\Z}$ is a mapping from $E$ to $E$ satisfying some suitable conditions, and $\sig(t)=(\sig_m(t))_{m\in\Z}$ is a function from $\R$ to $E$.

Let the function $\sig_0(t)\in L_c^2(\R,E)$ and take $$\Sig:=\Sig(\sig_0)=\mb{the closure of $\{\sig_0(s+\cdot):s\in \R\}$ in $L_{loc}^2(\R,E)$}$$
as the symbol space. Define the translation group $\theta$ on $\Sig$ as
$$\theta_t\sig=\sig(t+\cdot),\hs t\in \R,\, \sig\in\Sig.$$

Assume that for each fixed $u_0\in E, \sig\in \Sig$ and $t_0\in\R$, equations \eqref{e3.1}-\eqref{e3.2} has a unique continuous solution $u(t)=u(t,t_0;u_0,\sig)$ on $[t_0,+\8)$ in the sense that $$u:[t_0,+\8)\ra E\,\mb{is continuous with $u(t_0)=u_0$, and satisfies the equation \eqref{e3.1}.}$$ Moreover, for $t_0=0$, operators of solutions
$$U_\sig(t,0):u_0\mapsto u(t,0;u_0,\sig)=U_\sig(t,0)u_0,\hs t\geq0,\,u_0\in E,\,\sig\in \Sig.$$
generate a family of continuous processes $\{U_\sig(t,0)\}_{t\geq 0},\sig\in \Sig$ in $E$. Then one can easily check that the assumption ({\bf I}) in Section 2 holds true for $\{U_\sig(t,0)\}_{t\geq 0},\sig\in \Sig$.

\vs
For convenience in statement, we define a semigroup $\Psi=\{\Psi(t)\}_{t\geq0}$, called a skew-product flow on $E\X \Sig$ corresponding to $\{U_\sig(t,0)\}_{t\geq 0},\sig\in \Sig$ by
\be\label{e3.3}\Psi(t)(u,\sig)=(U_\sig(t,0)u,\theta_t\sig),\hs (u,\sig)\in E\X \Sig .\ee

In the following, we state and prove our main results in this section. We first construct a family of sets $\{\cA_\ve(\sig)\}_{\sig\in\Sig}$ in an arbitrary small $\ve$-neighborhood of the global attractor of
$\Psi$, such that $\{\cA_\ve(\sig)\}_{\sig\in\Sig}$ uniformly (w.r.t. $\sig\in \Sig$) forward attracts each bounded subset of $E$.

\vs
\bt\label{t3.1}
Assume that the family of processes $\{U_\sig(t,0)\}_{t\geq 0},\sig\in \Sig$ has a bounded uniformly (w.r.t. $\sig\in \Sig$) absorbing set $B_0$. Moreover,
$$\lim_{t\ra+\8}\lim_{M\ra+\8}\sup_{u\in B_0}\sum_{|m|>M}|\big(U_\sig(t,0)u\big)_m|^2=0,\Hs \forall \sig\in \Sig.$$
Then for any $\ve>0$, there exists a family of sets $\{\cA_\ve(\sig)\}_{\sig\in\Sig}$ with $$A(\sig)\subset\cA_\ve(\sig)\subset \cN_\ve(\sig),\Hs \A\sig\in \Sig,$$ such that it is forward invariant for $\{U_\sig(t,0)\}_{t\geq 0},\sig\in \Sig$, where $\cN_\ve(\sig)$ denotes the $\sig$-section of an $\ve$-neighborhood of $\Cup\limits_{\sig\in \Sig}(A(\sig)\X\{\sig\})$. Furthermore, $\{\cA_\ve(\sig)\}_{\sig\in\Sig}$ uniformly (w.r.t. $\sig\in \Sig$) forward attracts each bounded set $B\ss E$: \be\label{e3.4}\lim_{t\ra +\infty}\sup_{\sig\in \Sig} {\rm d}_H(U_\sig(t,0)B,\cA_\ve(\theta_t\sig))=0.\ee
\et
\bo
Since $B_0$ is a bounded uniformly (w.r.t. $\sig\in\Sig$) absorbing set for the family of processes $\{U_\sig(t,0)\}_{t\geqslant 0},\sig\in \Sig$ in $E$ and
$$\lim_{t\ra+\8}\lim_{M\ra+\8}\sup_{u\in B_0}\sum_{|m|>M}|\big(U_\sig(t,0)u\big)_m|^2=0,\Hs \forall \sig\in \Sig,$$
one can conclude from \cite[Theorem 3.1]{ZZL} that $\{U_\sig(t,0)\}_{t\geqslant 0},\sig\in\Sig$ has a compact uniform attractor $D=\w_\Sig(B_0)$ in $E$.

Now let $\Psi$ be the semigroup given by \eqref{e3.3}. Then we infer from \cite{CLR,Chep} that $\Psi$ has a global attractor $\cA$ on $ E \times\Sig$ with \be\label{e3.5}\cA=\bigcup_{\sig\in \Sig}\big(A(\sig)\times\{\sig\}\big),\ee where $\{A(\sig)\}_{\sig\in \Sig}$ is the pullback attractor of $\{U_\sig(t,0)\}_{t\geqslant 0},\sig\in\Sig$ and $D=\Cup_{\sig\in\Sig}A(\sig)$.

Let $\ve>0$ be fixed. Write $\cD= D\times \Sig$. Noticing that the global attractor $\cA$ attracts $\cD$, there exists $T:=T(\ve,\cD)>0$ independent of $\sig\in \Sig$ such that $$\Psi(t)\cD\subset N_\ve(\cA),\Hs t\geqslant T.$$
Set $\cA_\ve=\bigcup_{s\geqslant T}\Psi(s)\cD$. By the definition of semigroups, one can easily verify that the set $\cA_\ve$ is $\Psi$-positively invariant and the $\sig$-section $\cA_\ve(\sig)$ of $\cA_\ve$ satisfies
$$\cA_\ve(\sig)=\bigcup_{s\geqslant T}U_\sig(0,-s)D.$$ That is
$$\ba{l}\Psi(t)\Cup_{\sig\in \Sig}\big(\cA_\ve(\sig)\X\{\sig\}\big)\subset \Cup_{\sig\in \Sig}\big(\cA_\ve(\sig)\X\{\sig\}\big),\hs \A t\geqslant 0.\ea$$ Hence
$$\ba{l}\Cup_{\sig\in \Sig}\big(U_\sig(t,0)\cA_\ve(\sig)\X\{\theta_t \sig\}\big)\subset \Cup_{\sig\in \Sig}\big(\cA_\ve(\theta_t\sig)\X\{\theta_t\sig\}\big)\hs t\geqslant 0,\ea$$ which implies the family of sets $\{\cA_\ve(\sig)\}_{\sig\in \Sig}$ is forward invariant for $\{U_\sig(t,0)\}_{t\geqslant 0},\sig\in\Sig$. By the invariance of the pullback attractor $\{A(\sig)\}_{\sig\in \Sig}$, one has $A(\sigma) \subset {\mathcal A}_\varepsilon(\sig),\sig\in \Sig$. We infer from the construction of $\cA_\ve$ that $\cA_\ve\subset N_\ve(\cA)$, which shows $$\Cup_{\sig\in \Sig}\big(\cA_\ve(\sig)\times\{ \sig\}\big)\subset \Cup_{\sig\in \Sig}\big(\cN_\ve(\sig)\times\{\sig\}\big),$$ where $\cN_\ve(\sig)$ denotes the $\sig$-section of $N_\ve(\cA)$. Thereby ${\mathcal A}_\varepsilon(\sig)\subset {\mathcal N}_\varepsilon(\sig)$ for $\sig\in \Sig$.
Thus $$A(\sigma) \subset {\mathcal A}_\varepsilon(\sig)\subset {\mathcal N}_\varepsilon(\sig),\hs\sig\in\Sig.$$
 Write  $\cA_\ve^T(\sig)=U_\sig(0,-T)D$. Then the set-valued mapping $\sig\ra \cA_\ve^T(\sig)$ is continuous at $\sig\in \Sig$. Indeed, for each $u\in D$ and the fixed $T>0$, by the continuity of processes, we see that $U_\cdot(0,-T)u$ is continuous from $\Sig$ to $E$. Since the set $D$ is compact, one deduces that the mapping $U_\cdot(0,-T)D$ is continuous from $\Sig$ to $E$ in the sense of Hausdorff distance. Thereby it follows from the construction of $\mathcal{A}^T_\ve(\sig)$ that the mapping $\sig\rightarrow \cA_\ve^T(\sig)$ is continuous at $\sig\in \Sig$.
\vs
In the following we shall prove that the family of sets $\{\cA_\ve(\sig)\}_{\sig\in \Sig}$ uniformly forward attracts each bounded set $B\ss E$.
For any $\sig\in\Sigma$ and $t\geq T$, we observe that
\begin{align}\label{ee3.6}
{\rm d}_H(U_\sig(t,0)B,\cA_\varepsilon(\theta_t\sig))&\leq{\rm d}_H(U_\sig(t,0)B,\cA_\varepsilon^T(\theta_t\sig))\nonumber \\
&={\rm d}_H(U_\sig(t,0)B,U_\sig(t,t-T)D)\nonumber \\
&={\rm d}_H(U_{\theta_{t-T}\sig}(T,0)U_\sig(t-T,0)B,U_{\theta_{t-T}\sig}(T,0)D).
\end{align}
Recalling that $D$ is the uniform attractor, we have
\be\label{ee3.7}
\sup_{\sig\in \Sig}{\rm d}_H(U_\sig(t-T,0)B,D)\ra 0\hs \mb{as} \hs t\ra \8.\ee
Therefore, by the continuity of $U_\cdot(T,0)$ on $D$, we conclude from \eqref{ee3.6} and \eqref{ee3.7} that the result holds true.
\eo

\br\label{r3.3}
In the proof of the above theorem, one sees that for each $\ve>0$, $$A(\sig)\subset \cA_\ve(\sig)\subset \cN_\ve(\sig),\hs \sig\in\Sig,$$ which implies $$\de_H(\cA_\ve(\sig),A(\sig))\ra 0 \hs \mb{as}\hs \ve\ra 0.$$ Moreover, the family of sets $\{\cA_\ve(\sig)\}_{\sig\in\Sig}$ is forward attracting. Thus the sets $\{\cA_\ve(\sig)\}_{\sig\in\Sig}$ may be used to describe the forward dynamical behavior of nonautonomous LDSs.
\er
\br\label{r3.4}
The continuity of the mapping $\sig\ra \cA_\ve^T(\sig)$ shows some synchronous properties of the set $\cA_\ve^T(\sig)$ with $\sig$. For example, if $\sig\in \Sig$ is (almost) periodic for $t\in \mathbb{R}$, then the mapping $t\ra \cA_\ve^T(\theta_t\sig)$ is (almost) periodic as well.
\er

Next, we continue to construct a family of sets $\{\cB_\ve(\sig)\}_{\sig\in\Sig}$ such that it forward exponentially attracts bounded sets. For this purpose, we make the following assumptions:
\vs
\noindent{\bf (A)} 
There is a closed neighborhood $\{\cN(\sig)\}_{\sig\in \Sig}$ of the pullback attractor $\{A(\sig)\}_{\sig\in \Sig}$ such that
\benu
\item[(i)] There exists  $\tau^*>0$ such that for any $\sig\in \Sig$, the operator $U_\sig(\tau^*,0)$ is a compact perturbation of the contraction on $\cN(\sig)$:
   $$\|U_\sig(\tau^*,0)u-U_\sig(\tau^*,0)v\|_E\leqslant \de\|u-v\|_E+\|K(\sig)u-K(\sig)v\|_E$$ for any $u,v\in \cN(\sig)$, where $0<\de<1/2$ and $K(\sig)$ is an operator from $\cN(\sig)$ to $F$, which is a Banach space compactly embedded into $E$ and satisfies $$\|K(\sig)u-K(\sig)v\|_F\leqslant L_1\|u-v\|_E,\hs u,v\in \cN(\sig)$$ with $L_1>0$ independent of $\sig.$
\item[(ii)] There is a positive constant $L_2$ (independent of $\tau$ and $\sig$) satisfying $$\|U_\sig(\tau,0)u-U_\sig(\tau,0)v\|_E\leqslant L_2\|u-v\|_E,\hs u,v\in \cN(\sig)$$ for any $\tau\in[0,\tau^*]$ and $\sig\in \Sig$.
\eenu
 Then we have
\bt\label{t3.4}
Assume that the conditions in Theorem \ref{t3.1} are satisfied and let the assumption {\bf (A)} hold true. Then for each $\ve>0$, there exists a family of sets $\{\cB_\ve(\sig)\}_{\sig\in\Sig}$ such that $\{\cB_\ve(\sig)\}_{\sig\in\Sig}$ is forward invariant under the acting of $\{U_\sig(t,0)\}_{t\geq 0},\sig\in \Sig.$ Moreover, the system has a pullback exponential attractor $\{\sM(\sig)\}_{\sig\in \Sig}$ with $$\sM(\sig)\ss\cB_\ve(\sig),\hs \A \sig\in \Sig.$$ Furthermore, the family of sets $\{\cB_\ve(\sig)\}_{\sig\in\Sig}$ uniformly forward exponentially attracts each bounded set $B\ss E$.
\et
\bo
First, let $\cA$ be the global attractor of the semigroup $\Psi$ given by \eqref{e3.3}. Write $\cB_0=B_0\X \Sig$. Then for each fixed $\ve>0$, there exists $T:=T(\cB_0,\ve)>0$ such that $$\Psi(t)\cB\subset \cN_\ve(\cA),\Hs t\geqslant T.$$
Set $\cB_\ve=\bigcup_{s\geqslant T}\Psi(s)\cB_0$. Similar to $\{\cA_\ve(\sig)\}_{\sig\in \Sig}$, one can also construct a family of sets $\{\cB_\ve(\sig)\}_{\sig\in \Sig}$ defined by $$\cB_\ve(\sig)=\bigcup_{s\geqslant T}U_\sig(0,-s)B_0,$$ such that  $\{\cB_\ve(\sig)\}_{\sig\in \Sig}$ is forward invariant for $\{U_\sig(t,0)\}_{t\geq 0},\sig\in \Sig$ and $\cB_\ve(\sig)\ss \cN_\ve(\sig)$ for all $\sig\in\Sig$.

Noticing that $\Cup_{\sig\in \Sig}\(\cN(\sig)\X\{\sig\}\)$  is a closed neighborhood of the global attractor $\cA$, we see that for each bounded set $B\subset E$, there is $T_B>0$ independent of $\sig\in \Sig$ such that
$$\Psi(t)\cB\subset \Cup_{\sig\in \Sig}\(\cN(\sig)\X\{\sig\}\)\hs t\geqslant T_B,$$ where $\cB=B\X\Sig$. It then follows from the definition of $\Psi(t)$ that if $t\geqslant T_B$,
$$\Cup_{\sig\in \Sig}\(U_\sig(t,0)B\X\{\theta_t \sig\}\)\subset \Cup_{\sig\in \Sig}\(\cN(\theta_t\sig)\X\{\theta_t\sig\}\),$$ from which it can be deduced that the family of sets $\{\cN(\sig)\}_{\sig\in \Sig}$ is uniformly (w.r.t. $\sig\in\Sig$) absorbing for the family of processes $\{U_\sig(t,0)\}_{t\geq 0},\sig\in \Sig.$ Hence we may suppose that the $T$ in the construction of $\{\cB_\ve(\sig)\}_{\sig\in \Sig}$ is sufficiently large so that $$\cB_{\ve}(\sig)\subset \cN(\sig)\hs \mb{for all $\sig\in \Sig$}.$$
Since the set $B_0$ in the construction of $\cB_\ve(\sig)$ is uniformly absorbing, the family of sets $\{\cB_\ve(\sig)\}_{\sig\in \Sig}$ is uniformly absorbing as well. Clearly, $\{\cB_\ve(\sig)\}_{\sig\in \Sig}$ satisfies conditions (i) and (ii) in assumption $\mathbf{(A)}$ and is uniformly bounded with respect to $\sig\in \Sig$.

Now, we show that the family of processes $\{U_\sig(t,0)\}_{t\geq 0},\sig\in \Sig$ has a pullback exponential attractor $\{\sM(\sig)\}_{\sig\in \Sig}$. Note that $\theta_t$ is a homeomorphism for each $t\in\R$. Then for any $\sig\in\Sig$, there exists $s\in \R$ such that $\sig=\theta_s\sig_0$.
Define \be\label{e3.7} U(t,s)u=U_{\sig_0}(t,s)u,\hs t\geqslant s,\hs u\in E,\ee where $\sig_0$ is the function given in  Section 3.  It is easy to see that $\{U(t,s)\}_{t\geqslant s}$ is a process in $E$. Set \be\label{e3.8}\cB_\ve(t):=\cB_\ve(\theta_t\sig_0),\hs t\in\R.\ee Then $\{\cB_\ve(t)\}_{t\in\R}$ is positively invariant for the process $\{U(t,s)\}_{t\geqslant s}$.

In what follows we first check that the family of sets $\{\cB_\ve(s)\}_{s\in \R}$ satisfies the conditions (1)-(5) in \cite[p. 656]{EYY}, and hence by \cite[Theorem 2.1]{EYY} the process $\{U(t,s)\}_{t\geqslant s}$ has a pullback exponential attractor $\{\sM(t)\}_{t\in\R}$ with $\sM(t)\subset \cB_\ve(t)$ for all $t\in\R$.
Indeed, observing that $\{\cB_\ve(\sig)\}_{\sig\in \Sig}$ satisfies conditions (i) and (ii) in assumption $\mathbf{(A)}$, we conclude that for each $\sig\in \Sig$ with $\sig=\theta_s\sig_0$ for some $s\in \R$, and the $\tau^*>0$,
$$\|U_{\sig_0}(\tau^*+s,s)u-U_{\sig_0}(\tau^*+s,s)v\|_E\leqslant \de\|u-v\|_E+\|K(\theta_s\sig_0)u-K(\theta_s\sig_0)v\|_E$$ for any $u,v\in \cB_\ve(\theta_s\sig_0)$ and $$\|K(\theta_s\sig_0)u-K(\theta_s\sig_0)v\|_F\leqslant L_1\|u-v\|_E,\hs u,v\in \cB_\ve(\theta_s\sig_0).$$ Moreover, for any $\tau\in[0,\tau^*],$
$$\|U_{\sig_0}(\tau+s,s)u-U_{\sig_0}(\tau+s,s)v\|_E\leqslant L_2\|u-v\|_E,\hs u,v\in \cB_\ve(\theta_s\sig_0).$$
By \eqref{e3.7} and \eqref{e3.8}, one finds that
$$\|U(\tau^*+s,s)u-U(\tau^*+s,s)v\|_E\leqslant \de\|u-v\|_E+\|K(s)u-K(s)v\|_E$$ for any $u,v\in \cB_\ve(s)$,
$$\|K(s)u-K(s)v\|_F\leqslant L_1\|u-v\|_E,\hs u,v\in \cB_\ve(s)$$
and for any $\tau\in[0,\tau^*]$,
$$\|U(\tau+s,s)u-U(\tau+s,s)v\|_E\leqslant L_2\|u-v\|_E,\hs u,v\in \cB_\ve(s),$$ which shows that $\{\cB_\ve(s)\}_{s\in \R}$ satisfies the conditions (4)-(5) in \cite[p. 656]{EYY} for the process $\{U(t,s)\}_{t\geqslant s}$ defined by \eqref{e3.7}. Similarly, it is easy to verify that $\{\cB(s)\}_{s\in \R}$ satisfies conditions (1)-(3) in \cite[p. 656]{EYY}.
 Therefore one can immediately conclude from \cite[Theorem 2.1]{EYY} that the process $\{U(t,s)\}_{t\geqslant s}$ has a pullback exponential attractor $\{\cM(t)\}_{t\in\R}$ with $\cM(t)\subset \cB_\ve(t)$ for all $t\in \R$, and \be\label{e3.9}{\rm d}_H\big(U(t,s)\cB_\ve(s),\cM(t)\big)\leqslant C{\rm e}^{-\a(t-s)},\hs s\in \R,\,s\leqslant t<+\8,\ee
where $C,\a$ are some positive constants independent of $s\in\R$, and $C$ depends on the boundedness of the family of sets $\{\cB_\ve(s)\}_{s\in \R}$ uniformly for $s\in\R$.

For any $s\in \R$, define\be\label{e3.10}\sM(\sig)=\cM(s),\hs \sig=\theta_s\sig_0.\ee
Clearly, $\sM(\sig)\subset \cB_\ve(\sig),\sig\in \Sig.$
We shall show that $\{\sM(\sig)\}_{\sig\in \Sig}$ is a pullback exponential attractor of $\{U_\sig(t,0)\}_{t\geq 0},\sig\in \Sig$. First, since $\theta_t:\Sig\mapsto \Sig$ is a homeomorphism for $t\in \R$, for any $\sig\in \Sig$, there exists $s\in \R$ such that $\sig=\theta_{s}\sig_0$. Thus $\sM(\sig)$ is compact, as the set $\cM(s)$ is compact. For the fractal dimension of $\sM(\sig)$, noticing that the fractal dimension of $\cM(t)$ is finite and uniformly bounded for $t\in \R$, one can conclude that the fractal dimension of $\cM(\sig)$ is finite and uniformly bounded with respect to $\sig\in \Sig$ as well. Now we verify that the family of sets $\{\sM(\sig)\}_{\sig\in \Sig}$ is forward invariant for $\{U_\sig(t,0)\}_{t\geq 0},\sig\in \Sig$. Let $\sig=\theta_s\sig_0$ for some $s\in\R$ and $\tau\geqslant 0$. We infer from \eqref{e3.7}, \eqref{e3.10} and the positive invariance of $\{\cM(t)\}_{t\in \R}$ for $\{U(t,s)\}_{t\geqslant s}$ that
\begin{align*}
U_\sig(\tau,0)\sM(\sig)=U_{\sig_0}(\tau+s,s)\sM(\theta_s\sig_0)&=U(\tau+s,s)\cM(s)\\
&\subset \cM(\tau+s)=\sM(\theta_{\tau+s}\sig_0)=\sM(\theta_\tau \sig).
\end{align*}
Finally, it remains to prove the exponential attractivity of $\{\sM(\sig)\}_{\sig\in \Sig}.$ For each $B\subset E$, we conclude from the definition of the pullback exponential attractor $\{\cM(t)\}_{t\in \R}$ for $\{U(t,s)\}_{t\geqslant s}$ that if $\tau\geqslant 0$ and $\sig=\theta_s\sig_0$, there exist constants $\alpha>0$, $C_B$ and $T_{\sig,B}$ depending on $T_{s,B}$ such that
\begin{align*}
{\rm d}_H\(U_{\sig}(0,-\tau)B,\sM(\sig)\)&={\rm d}_H(U_{\sig_0}(s,s-\tau)B,\sM(\theta_{s}\sig_0))\\
&={\rm d}_H(U(s,s-\tau)B,\cM(s))
\leqslant C_B {\rm e}^{-\alpha \tau}, \hs \tau\geqslant T_{\sig,B}.
\end{align*}
Thus, the family of sets $\{\sM(\sig)\}_{\sig\in \Sig}$ is a pullback exponential attractor for the family of processes $\{U_\sig(t,0)\}_{t\geq 0},\sig\in \Sig$ with $\sM(\sig)\ss\cB_\ve(\sig)$ for all $\sig\in\Sig$.

To complete the proof of the theorem, it remains to prove the forward exponential attractivity of $\{\cB_\ve(\sig)\}_{\sig\in\Sig}$.
By the construction of a pullback exponential attractor (see \eqref{e3.9}) and  \eqref{e3.7}, \eqref{e3.8}, \eqref{e3.10} that
\be\label{e3.11}
{\rm d}_H\(U_\sig(0,-t)\cB_\ve(\theta_{-t} \sig),\sM(\sig)\)\leqslant C{\rm e}^{-\a t},\Hs \A \sig\in \Sig,\,\, \A t\geqslant 0.\ee
Let $B\subset E$ be a bounded set. Since $\{\cB_\ve(\sig)\}_{\sig\in \Sig}$ is uniformly absorbing, it follows from \eqref{e3.11} that there exists $T_B>0$ (independent of $\sig$) such that
$${\rm d}_H\(U_\sig(0,-t)B,\sM(\sig)\)\leqslant C_1{\rm e}^{-\a t},\hs \A t\geqslant T_B,$$ where $C_1$ is a positive constant depending on $C$ and $T_B$. Recalling that $\sM(\sig)\subset \cB_\ve(\sig)$ for all $\sig\in \Sig$, one can deduce that $${\rm d}_H\(U_\sig(0,-t)B,\cB_\ve(\sig)\)\leqslant {\rm d}_H\(U_\sig(0,-t)B,\sM(\sig)\).$$ Hence we conclude that if $t\geqslant T_B$, then $${\rm d}_H\(U_\sig(0,-t)B,\cB_\ve(\sig)\)\leqslant C_1{\rm e}^{-\a t},\hs \A \sig\in\Sig,$$
from which it can be seen that
$$\sup_{\sig\in \Sig}{\rm d}_H\(U_\sig(0,-t)B,\cB_\ve(\sig)\)\leqslant C_1{\rm e}^{-\a t},\hs \A t\geqslant T_B.$$
Therefore $$\sup_{\sig\in \Sig}{\rm d}_H\(U_\sig(t,0)B,\cB_\ve(\theta_t \sig)\)\leqslant C_1{\rm e}^{-\a t},\hs\A t\geqslant T_B.$$
The proof of the theorem is complete.
\eo

\br\label{r3.5}
In the proof of the above theorem, one can see that $$ A(\sig)\subset \sM(\sig)\subset \cB_\ve(\sig)\subset \cN_\ve(\sig),\hs \A \sig\in \Sig,$$  which implies
$$\mathrm{\de}_H(\sM(\sig),A(\sig))\ra 0,\hs \de_H(\cB_\ve(\sig),\sM(\sig))\ra 0 \hs \mb{as} \hs \ve\ra 0.$$
\er

We infer from the proof of Theorem \ref{t3.4} that the existence of a pullback exponential attractor plays an important role in proving the exponential attraction of $\{\cB_\ve(\sig)\}_{\sig\in\Sig}$. Therefore by using the results on the existence of pullback exponential attractors in a Hilbert space, see \cite[Theorem 2]{ZH}, we can obtain the following main results, which seem to be more convenient in applications.
\bt\label{t3.6}
Let $\{U_\sig(t,0)\}_{t\geq 0},\sig\in \Sig$ be the family of processes generated by equations \eqref{e3.1}-\eqref{e3.2} in $H$, where $H$ is a Hilbert space. Assume that $B_0\ss H$ is a bounded uniformly (w.r.t. $\sig\in\Sig$) absorbing set for $\{U_\sig(t,0)\}_{t\geq 0},\sig\in \Sig$ and satisfies
$$\lim_{t\ra+\8}\lim_{M\ra+\8}\sup_{u\in B_0}\sum_{|m|>M}|\big(U_\sig(t,0)u\big)_m|^2=0,\Hs \forall \sig\in \Sig,$$
and for some $T_0:=T(B_0)>0$ such that $$U_\sig(t,0)B_0\subset B_0,\hs \A t\geqslant T_0,\,\, \A \sig\in \Sig.$$
In addition, if
\benu
\item[{\rm ({\bf A}1)}] there exist positive constants $T^*>T_0$ and $L=L_{T^*}$ such that for each $\sig\in \Sig$, $$\|U_\sig(t,0)u-U_\sig(t,0)v\|_H\leqslant L\|u-v\|_H,\hs u,v\in B_0,\,\,t\in [T_0,T^*];$$
\item[{\rm ({\bf A}2)}] there are $0\leqslant \beta<1/2$ and a finite dimensional subspace $H_m$ of $H$ such that for any $\sig\in \Sig$, $$\|(I-P_m)(U_\sig(T^*,0)u-U_\sig(T^*,0)v)\|_H\leqslant \beta\|u-v\|_H,\hs u,v\in B_0,$$ where $P_m:H\mapsto H_m$ is a projection, $\beta$ and $m\in\N$ (depending on $T^*$) are independent of $\sig\in \Sig$,
\eenu
then for each $\ve>0$, there exists a family of sets $\{\cB_\ve(\sig)\}_{\sig\in \Sig}$ such that  $\{\cB_\ve(\sig)\}_{\sig\in \Sig}$ is forward invariant for $\{U_\sig(t,0)\}_{t\geq 0},\sig\in \Sig$ and $\{\cB_\ve(\sig)\}_{\sig\in \Sig}$ uniformly forward exponentially attracts each bounded subset of $H$.
\et
\bo
For each fixed $\ve>0$, let $\{\cB_\ve(\sig)\}_{\sig\in\Sig}$ be the family of sets given in the proof of Theorem \ref{t3.4}, where $$\cB_\ve(\sig)=\bigcup_{s\geqslant T}U_\sig(0,-s)B_0.$$
Since the set $B_0$ is uniformly absorbing, it may be assumed that the $T$ in the construction of $\cB_\ve(\sig)$ is sufficiently large so that $$U_\sig(0,-s)B\ss B_0,\hs s\geqslant T$$ uniformly for $\sig\in \Sig$. Thus $\{\cB_\ve(\sig)\}_{\sig\in\Sig}$ satisfies the conditions ({\bf A}1)-({\bf A}2).

By \eqref{e3.7}, we infer from conditions ({\bf A}1)-({\bf A}2) that for $\sig=\theta_s\sig_0$, $$\|U(t+s,s)u-U(t+s,s)v\|_H\leqslant L\|u-v\|_H,\hs u,v\in B_0,\,\,t\in [T_0,T^*]$$ and $$\|(I-P_m)(U(T^*+s,s)u-U(T^*+s,s)v)\|_H\leqslant \beta\|u-v\|_H,\hs u,v\in B_0.$$
Hence by \cite[Theorem 2]{ZH}, one can conclude that the process $\{U(t,s)\}_{t\geqslant s}$ has a pullback exponential attractor $\{\cM(t)\}_{t\in\R}$.
Note that the family of sets $\{\cB_\ve(\sig)\}_{\sig\in\Sig}$ is forward invariant and uniformly absorbing for $\{U_\sig(t,0)\}_{t\geq 0},\sig\in \Sig$. Repeating the same argument below \eqref{e3.9}, it can be shown that $\{\cB_\ve(\sig)\}_{\sig\in\Sig}$ uniformly forward exponentially attracts each bounded subset of $H$.
\eo

\setcounter {equation}{0}

\section{Applications}

In this section, we study the nonautonomous discrete Gray-Scott equations driven by quasi-periodic external forces, illustrating how to apply our abstract results to some concrete equations.

Consider the following nonautonomous LDS:
 \begin{align}
\dot{u}_{m}
&=-d_1(Au)_m-(\lambda+k)u_m+u_m^2v_m-\alpha
u_m^3+\beta z_m+b_{1m}f_{1m}(\sigma(t)),\label{aa1}\\
  \dot{v}_{m}
&=-d_2(Av)_m+\lambda(a_m-v_m)-u_m^2v_m+\alpha
u_m^3+b_{2m}f_{2m}(\sigma(t)), \label{aa2}\\
  \dot{z}_{m}
&=-d_3(Az)_m+ku_m-(\lambda+\beta)z_m+b_{3m}f_{3m}(\sigma(t))\label{aa3}
 \end{align}
for $m\in \mathbb{Z},t>t_0, t_0\in \mathbb{R}$,
associated with the initial conditions
\begin{eqnarray}\label{aa4}
 u_{m}(t_0)=u_{m,0},\hs v_{m}(t_0)=v_{m,0},\hs z_{m}(t_0)=z_{m,0},\Hs m\in
 \mathbb{Z},
\end{eqnarray}
where $d_1,d_2,d_3,k,\alpha,\beta,\lambda$ are positive numbers,
$a_m\in \mathbb{R},b_{im}\in\mathbb{R},f_{im}\in
\mathbb{R},i=1,2,3,\sigma\in \mathbb{T}^\kappa(\kappa$-dimensional
torus), and $A$ is a linear operator defined by
\be\label{aa5}(Au)_{m}=2u_{m}-u_{m+1}-u_{m-1},\Hs \A u=(u_m)_{m\in\Z}.\ee
Equations \eqref{aa1}-\eqref{aa3} can be regarded as a discrete analogue of the
following nonautonomous three-component reversible Gray-Scott model on
$\mathbb{R}$:
\begin{align}
  \frac{\partial u}{\partial t}
=&
  d_1u_{xx}-(\lambda+k)u+u^2v-\alpha u^3+\beta z+b_1f_1(x,\sigma),\label{aa6}\\
  \frac{\partial v}{\partial t}
=&
  d_2v_{xx}- \lambda v-u^2v+\alpha u^3+\lambda a+b_2f_2(x,\sigma),\label{aa7}\\
  \frac{\partial z}{\partial t}
=&
  d_3z_{xx}+ku-(\lambda+\beta)z+b_3f_3(x,\sigma),\label{aa8}
\end{align}
where the positive numbers $d_1,d_2,d_3,k,\alpha,\beta,\lambda$ are defined by nondimensionalization, see \cite{Y2} for details.

In order to express equations \eqref{aa1}-\eqref{aa4} as an abstract first-order ODE, we set $u=(u_m)_{m\in\Z}, v=(v_m)_{m\in\Z}, z=(z_m)_{m\in\Z}, u^2v=(u_m^2v_m)_{m \in \mathbb{Z}}, u^3=(u_m^3)_{m\in\Z},$
$b_if_i(\sigma(t))=(b_{im}f_{im}(\sigma(t)))_{m\in \Z}, i=1,2,3, a=(a_m)_{m\in
\Z},\sigma(t)=({\bf x}t+\sigma_0) \,{\rm mod}(\mathbb{T}^\kappa).$ Then
\begin{align}
 \dot{u}
 &=-d_1 Au-(\lambda+k)u+u^2v-\alpha u^3+\beta z+b_1f_1(\sigma(t)), \label{4.1} \\
\dot{v}
 &=-d_2 Av-\lambda v-u^2v+\alpha u^3+\lambda a+b_2f_2(\sigma(t)),\label{4.2}\\
 \dot{z}
 &=-d_3 Az+ku-(\lambda+\beta)z+b_3f_3(\sigma(t)),\label{4.3}
\end{align}
associated with initial conditions
\be\label{4.4}
u(t_0)=(u_{m,0})_{m\in \Z},\, v(t_0)=(v_{m,0})_{m\in \Z},\, z(t_0) =(z_{m,0})_{m\in \Z},\hs t_0\in\R,
\ee
where
 $\sigma_0\in \mathbb{T}^\kappa$ ({\bf x} and $\mathbb{T}^\kappa$ will be introduced below), and $A$ is a linear operator from $\ell^2$ to $\ell^2$ defined by \eqref{aa5}.


\subsection{Mathematical setting}
We first introduce some operators and spaces.

Let $\ell^2$ be the space given in Section 3. Define the linear operators $B$ and $B^*$ on $\ell^2$ as
$$
  (Bu)_{m}=u_{m+1}-u_{m},\,\, (B^{*}u)_{m}=u_{m-1}-u_{m},
   \hs  m\in \mathbb{Z},\,\A u=(u_m)_{m\in\Z}\in \ell^2.\label{2.3}
$$
Then it is easy to see that $B^{*}$ is the adjoint operator of $B$. Moreover, for any $u,v \in \ell^2$,
\begin{align*}
  (Au,v)&=(B^{*}Bu,v)=(Bu,Bv), \hs
  (Bu,v)=(u,B^{*}v), \\
  \|Bu\|^{2}&=\|B^{*}u\|^{2}\leqslant 4\|u\|^{2},\hs
  \|Au\|^{2}\leqslant 16\|u\|^{2}.
\end{align*}

Let $\mathbb{T}^\kappa$ be a $\kappa$-dimensional torus
$$
  \mathbb{T}^\kappa=\{\sigma=(\sigma_1,\cdot\cdot\cdot,\sigma_\kappa):\sigma_i\in [-\pi,\pi],\hs\forall i=1\cdot\cdot\cdot \kappa\}
  $$
with the following identification
$$
  (\sigma_1,\cdot\cdot\cdot,\sigma_i,-\pi,\sigma_{i+1},\cdot\cdot\cdot\sigma_\kappa)\sim
  (\sigma_1,\cdot\cdot\cdot,\sigma_i,\pi,\sigma_{i+1},\cdot\cdot\cdot\sigma_\kappa),\Hs
  \forall i=1,\cdot\cdot\cdot,\kappa,
$$
and the topology and metric induced from the topology and metric on
$\mathbb{R}^\kappa$. Thus we consider the norm on
$\mathbb{T}^\kappa$ given by
$$
  \|\sigma\|_{\mathbb{T}^\kappa}=(\sum\limits_{i=1}^{\kappa}\sigma_i^2)^{1/2},\Hs
  \forall \sigma=(\sigma_1,\cdot\cdot\cdot,\sigma_\kappa)\in \mathbb{T}^\kappa.
$$
Let $\mathbf{x}=(x_1,\cdot\cdot\cdot,x_\kappa)\in \mathbb{R}^\kappa$
be a fixed vector and $x_1,\cdot\cdot\cdot,x_\kappa$ be rationally
independent. For $s\in \mathbb{R}$, define
$$
  \theta_s\sigma_0=(\mathbf{x}s+\sigma_0){\rm mod}(\mathbb{T}^\kappa),\Hs \sigma_0\in
  \mathbb{T}^\kappa.
  $$
Then one can easily check that
$$
  \theta_s\mathbb{T}^\kappa=\mathbb{T}^\kappa,\Hs \forall s\in\mathbb{R}.
  $$
Denote by $E$ the product space $\ell^{2}\times \ell^{2}\times \ell^{2}.$
Then $E$ is a Hilbert space as well. We also use the notations
$(\.,\.)$ and $||\.||$ to denote the inner product and norm,
respectively. Define the extended space $\tilde{E}=E\times
\mathbb{T}^\kappa$ with the norm defined by
$$
  \|\~\vp\|=(\|\vp\|^2+\|\sigma\|^2_{\mathbb{T}^\kappa})^{1/2},\Hs \forall
  \~\vp=(\vp,\sigma)\in \tilde{E}.
  $$
Now we rewrite equations \eqref{4.1}-\eqref{4.4} as an abstract first-order ODE:
\begin{align}
&\dot{\varphi}+\Theta\varphi=G(\varphi,\sigma,t),\Hs t>t_0,\label{4.5}\\
 & \varphi(t_0)=\varphi_{0}=(u_{0},v_0,z_0)^{T}, \Hs t_0\in \mathbb{R},\label{4.6}
\end{align}
where $\varphi=(u,v,z)^{T}$, and
$$
  G(\varphi,\sigma,t)=(u^2v-\alpha
  u^3+b_1f_1(\sigma(t)),\lambda a-u^2v+\alpha u^3+b_2f_2(\sigma(t)),b_3f_3(\sigma(t)))^{T},$$
$$
\begin{array}{ccc}
  \Theta=\left(\begin{array}{ccc}
 d_1 A+(\lambda+k)I& 0 & -\beta I  \\ \noalign{\medskip}
 0  &  d_2 A+\lambda I & 0
 \\ \noalign{\medskip}
 -k I & 0  &  d_3 A+(\lambda+\beta)I
\end{array}\right).
\end{array}
$$
For the parameters $a,b_i$ and functions $f_i$ for $i=1,2,3$ in equations \eqref{4.1}-\eqref{4.3}, we make the following assumptions:
\vs
(\textbf{H}1) Assume $a=(a_m)_{m\in \Z}\in
\ell^2,b_i=(b_{i,m})_{m\in\mathbb{Z}}\in \ell^2, i=1,2,3$, and that
$$
  2\beta<k<\min\big\{\frac{3\lambda}{2\mu-1},\beta+\lam\},
$$
where $\mu=k/\beta$.
\vs

(\textbf{H}2) $f_{im}(0_{\mathbb{T}^\kappa})=0$ for all
$m\in\mathbb{Z}, i=1,2,3$, and there exist $c_1,c_2,c_3>0$ such that
$$
  |f_{im}(\sigma_{01})-f_{im}(\sigma_{02})|\leq
  c_i\|\sigma_{01}-\sigma_{02}\|_{\mathbb{T}^\kappa}, \Hs \forall m\in \mathbb{Z}, \,\sigma_{01},\sigma_{02}\in \mathbb{T}^\kappa, \hs i=1,2,3.
$$

\subsection{Unique existence and boundedness of solutions}
We prove the existence and uniqueness of solutions and
then verify the existence of bounded absorbing set for skew-product flows. For simplicity, we let
$$
  Z(t)=\frac{\beta}{k}z(t).
$$ Then equations \eqref{4.1}-\eqref{4.3} can be transformed into the following equivalent form
\begin{align}
 \dot{u}
 &=-d_1 Au-(\lambda+k)u+u^2v-\alpha u^3+k Z+b_1f_1(\sigma(t)), \label{4.7}\\
\dot{v}
 &=-d_2 Av-\lambda v-u^2v+\alpha u^3+\lambda a+b_2f_2(\sigma(t)),\label{4.8}\\
 \mu\dot{Z}
 &=-d_3\mu AZ+ku-(\mu\lambda+k)Z+b_3f_3(\sigma(t)).\label{4.9}
\end{align}
\bt\label{t4.1}
Let the assumptions $(\mathbf{H}1)$-$(\mathbf{H}2)$ hold. Then for
any initial value $\varphi_0=(u_0,v_0,z_0)\in E, t_0\in\R$  and
$\sigma\in \mathbb{T}^\kappa$, equations \eqref{4.5}-\eqref{4.6} has a unique
solution $\vp(t)=\varphi(t,t_0;\vp_0,\sig)$ with
$$
\vp(t)=(u(t),v(t),z(t))\in {\mathcal C}([t_0, +\infty),E)\cap {\mathcal C}^{1}((t_0,
  +\infty),E),\,\, t\geqslant t_0.
  $$
\et
\begin{proof}
Let
$$
  \tilde{G}(\varphi,\sigma,t)=G(\varphi,\sigma,t)-\Theta\varphi,\Hs
  \sigma\in \mathbb{T}^\kappa.
$$
Then under the assumptions, one can check that $\tilde{G}:
\tilde{E}\times \mathbb{R}\mapsto E$ is locally Lipschitz
continuous respect to $\vp$ for any $\sigma\in \mathbb{T}^\kappa$
and $t$ in any compact interval of $\mathbb{R}$. By the classical
theory of ODEs, we can obtain the existence and uniqueness of a
local solution of equations \eqref{4.5}-\eqref{4.6} and
$$
  \varphi(t)=(u(t),v(t),z(t))\in {\mathcal C}([t_0, T),E)\cap
{\mathcal C}^{1}((t_0,
  T),E)$$
for some $T>t_0$. Moreover, if $T<+\infty,$ then
$\lim\limits_{t\rightarrow T^-}\|\vp(t)\|=+\infty.$

In what follows, we prove that $T=+\infty.$ Taking the inner product
$(\cdot,\cdot)$ of the three equations of \eqref{4.7}-\eqref{4.9} with $\alpha
u(t),v(t),\mu\alpha Z(t)$, respectively, and summing up, we deduce
that
\begin{align}\label{4.10}
&
  \frac{1}{2}\frac{\mathrm{d}}{\mathrm{d}t}(\alpha\|u\|^2+\|v\|^2+\mu^2\alpha\|Z\|^2)
 +
  \alpha(\lam+k)\|u\|^2+\lam\|v\|^2+\mu\alpha(\mu\lam+k)\|Z\|^2\nonumber \\
\leqslant&\,
  k\alpha(1+\mu)\sum\limits_{m\in
  \mathbb{Z}}u_mZ_m+\lam\sum\limits_{m\in\mathbb{Z}}a_mv_m
 +\alpha\sum\limits_{m\in\mathbb{Z}}b_{1m}f_{1m}(\sigma)u_m+\sum\limits_{m\in\mathbb{Z}}b_{2m}f_{2m}(\sigma)v_m\nonumber \\
 &
 +\mu\alpha\sum\limits_{m\in\mathbb{Z}}b_{3m}f_{3m}(\sigma)Z_m.
\end{align}
Since
\begin{align}
  \lam\sum\limits_{m\in\mathbb{Z}}a_mv_m
&\leqslant
  \frac{\lam}{2}(\|a\|^2+\|v\|^2),\label{4.11}\\
  k\alpha(1+\mu)\sum\limits_{m\in \mathbb{Z}}u_mZ_m
&\leqslant
  \frac{k\alpha(1+\mu)}{2}(\|u\|^2+\|Z\|^2),\label{4.12}
  \end{align}
 and by $(\mathbf{H}1)$-$(\mathbf{H}2)$, we have
\begin{align}
  \alpha\sum\limits_{m\in\mathbb{Z}}b_{1m}f_{1m}(\sigma)u_m
&\leqslant
  \frac{\alpha(\lambda+k)}{4}\|u\|^2+\frac{\alpha
  \kappa c_1^2\pi^2}{\lam+k}\|b_1\|^2,\label{4.13}\\
  \sum\limits_{m\in\mathbb{Z}}b_{2m}f_{2m}(\sigma)v_m
&\leqslant
  \frac{\lambda}{4}\|v\|^2+\frac{\kappa c_2^2\pi^2}{\lam}\|b_2\|^2,\label{4.14}\\
  \mu\alpha\sum\limits_{m\in\mathbb{Z}}b_{3m}f_{3m}(\sigma)Z_m
&\leqslant
  \frac{\mu\alpha(\mu\lambda+k)}{4}\|Z\|^2+\frac{\mu\alpha
  \kappa c_3^2\pi^2}{\mu\lam+k}\|b_3\|^2.\label{4.15}
\end{align}
Thus it follows from \eqref{4.10}-\eqref{4.15} that for $t\geqslant t_0$,
\begin{align}\label{4.16}
  \frac{\mathrm{d}}{\mathrm{d}t}(\alpha\|u\|^2+\|v\|^2+\mu^2\alpha\|Z\|^2)
+
  \theta_1(\alpha\|u\|^2+\|v\|^2+\mu^2\alpha\|Z\|^2)
\leqslant
  R_0^2,
 \end{align}
 where
\begin{align*}
& \theta_1
 =
  \min\{\frac{3}{2}\lam-(\mu-\frac{1}{2})k,\frac{\lam}{2},\frac{3}{2}\lam+\frac{(\frac{\mu}{2}-1)k}{\mu^2}\},\\
& R_0^2 =
  \lam\|a\|^2+\frac{2\alpha \kappa c_1^2\pi^2}{\lam+k}\|b_1\|^2+\frac{2\kappa c_2^2\pi^2}{\lam}\|b_2\|^2+\frac{2\mu\alpha
  \kappa c_3^2\pi^2}{\mu\lam+k}\|b_3\|^2.
\end{align*}
Applying Gronwall inequality to \eqref{4.16} on $[t_0, t]$ with
$t\geqslant t_0$, we obtain \begin{align}\label{4.17}
\alpha\|u\|^2\!+\!\|v\|^2\!+\!\alpha\|z\|^2\leqslant \mathrm{e}^{-\theta_1
(t-t_0)}(\alpha\|u(t_0)\|^2\!+\!\|v(t_0)\|^2\!+\!\alpha\|z(t_0)\|^2)\!+\!\frac{R_0^2}{\theta_1}.
\end{align}
Let $\de_1=\min\{1,\alpha\}, \de_2=\max\{1,\alpha\}$. Then \eqref{4.17} shows that
$$
  \|\varphi(t)\|^2
\leqslant
  \frac{\delta_2}{\delta_1}\mathrm{e}^{-\theta_1
  (t-t_0)}\|\varphi(t_0)\|^2+\frac{R_0^2}{\theta_1\delta_1}. \Hs
  t\geqslant t_0,
$$
which implies that $T=+\infty$. The proof is complete.
\eo
Set $$U_\sig(t,0)\vp_0=\vp(t,0;\vp_0,\sig),\hs t\geqslant 0,\,\, \sig\in \T^\kappa.$$
Then one can easily verify that $\{U_\sig(t,0)\}_{t\geq 0},\sig\in \T^\kappa$ is a family of continuous processes in $E$.
Moreover, $\Psi$ is the corresponding semigroup. As a direct consequence of Theorem \ref{t4.1}, we can obtain the following
result.
\bl\label{l4.2} Suppose the assumptions $(\mathbf{H}1)$-$(\mathbf{H}2)$
hold. Then the family of processes $\{U_\sig(t,0)\}_{t\geq 0},\sig\in \T^\kappa$ generated by \eqref{4.5}-\eqref{4.6}
possesses a bounded uniformly $(w.r.t. \,\sig\in\T^\kappa)$ absorbing set $
B_0:=B(0,R)$ in $E$ centered at $0$ with radius
$\displaystyle R=2R_0/\sqrt{\theta_1\de_1}$.
Specifically, for any bounded set $B$ of $E$, there exists
$t(B)>0$ (independent of $\sig\in \T^\kappa$) such that
$$
  U_\sig(t,0)B \subset  B_0, \Hs\forall t\geqslant
  t(B),\hs \forall \sigma\in \mathbb{T}^\kappa.
$$
Particularly, there exists $T_0=t(B_0)>0$ such that
$$
  U_\sigma(t,0)B_0\subset B_0,\Hs\A t\geqslant T_0,\hs \A \sigma\in \mathbb{T}^\kappa.
  $$
\el
Next we prove that the solutions of \eqref{4.5}-\eqref{4.6} have uniformly
asymptotic end tails, from which one can also obtain the existence of a
uniform attractor.

\bl\label{l4.3}
Suppose the assumptions $(\mathbf{H}1)$-$(\mathbf{H}2)$ hold. Then
for any $\varepsilon>0$, there exist
$t_1(\varepsilon,B_0)>0$ and $N_1(\varepsilon,B_0)\in \N$,
such that the solution
$\varphi(t)=U_\sig(t,0)\vp_0=(u_m(t),v_m(t),z_m(t))_{m\in
\mathbb{Z}}$ $\in E$ of equations \eqref{4.5}-\eqref{4.6} with initial data
$\varphi_0 \in B_0$ satisfies
$$
  \sup\limits_{\varphi_0 \in B_0}\sum\limits_{|m|\geqslant N_1(\varepsilon,B_0)}|\varphi_{m}(t)|^{2}_{E}
\leqslant
  \varepsilon^2, \quad \forall t \geqslant t_1(\varepsilon,B_0),\hs \A\sigma\in \mathbb{T}^\kappa,
$$
where $|\varphi_{m}(t)|^{2}_{E}=u_m^2(t)+v_m^2(t)+z_m^2(t)$.
\el
\bo

Define a function $\chi(x)\in {\mathcal C}^{1}({\bf
\mathbb{R}}_+,[0,1])$ such that
\begin{eqnarray}
\quad \chi(x)=\left\{
\begin{array}{ll}\label{4.18}
    0, \enskip 0\leqslant  x\leqslant 1; \\
    1, \enskip x\geqslant 2,
\end{array}
\right. {\rm and }\, |\chi'(x)|\leqslant \chi_0 \, ({\rm positive
\enskip constant}), \, \forall x\geqslant 0.
\end{eqnarray}
Set
\begin{align*}
p=(p_m)_{m\in\mathbb{Z}}, \hs q=(q_m)_{m\in\mathbb{Z}},\hs
W=(W_m)_{m\in\mathbb{Z}},
\end{align*}
with
$$
  p_m=\chi(\frac{|m|}{M})u_m,\hs q_m=\chi(\frac{|m|}{M})v_m,\hs
  W_m=\chi(\frac{|m|}{M})Z_m,
$$
where $M$ is a positive integer. Taking the inner product
$(\cdot,\cdot)$ of these three equations of \eqref{4.7}-\eqref{4.9} with $\alpha
p, q$ and $\mu\alpha W,$ respectively, and adding these equations,
we obtain
\begin{align}\label{4.19}
& \frac{1}{2}\frac{\mathrm{d}}{\mathrm{d}t}\sum\limits_{m\in
  \mathbb{Z}}\chi(\frac{|m|}{M})(\alpha u_m^2\!+\!v_m^2\!+\!\mu^2\alpha
  Z_m^2)
 \!+\!
  d_1\alpha(Bu,Bp)\!+\!d_2(Bv,Bq)\!+\!\mu^2\alpha d_3(BZ,BW)\nonumber\\
&+\alpha(\lam+k)\sum\limits_{m\in
  \mathbb{Z}}\chi(\frac{|m|}{M})u_m^2
 +
  \lam\sum\limits_{m\in\mathbb{Z}}\chi(\frac{|m|}{M})v_m^2
 +
 \mu\alpha(\mu\lam+k)\sum\limits_{m\in
  \mathbb{Z}}\chi(\frac{|m|}{M})Z_m^2\nonumber\\
\leqslant&\,
  k\alpha(1+\mu)\sum\limits_{m\in
  \mathbb{Z}}\chi(\frac{|m|}{M})u_mZ_m+\lam\sum\limits_{m\in
  \mathbb{Z}}\chi(\frac{|m|}{M})a_mv_m+\alpha\sum\limits_{m\in
  \mathbb{Z}}\chi(\frac{|m|}{M})u_mb_{1m}f_{1m}(\sigma)\nonumber\\
&+\sum\limits_{m\in
  \mathbb{Z}}\chi(\frac{|m|}{M})v_mb_{2m}f_{2m}(\sigma)
 +
 \mu\alpha\sum\limits_{m\in
  \mathbb{Z}}\chi(\frac{|m|}{M})Z_mb_{3m}f_{3m}(\sigma).
\end{align}
Elementary computations show that
\begin{align}\label{4.20}
 (Bu,Bp)
&=
 \sum\limits_{m\in \mathbb{Z}}(Bu)_m(Bp)_m
 =\sum\limits_{m\in
 \mathbb{Z}}(Bu)_m[\chi(\frac{|m+1|}{M})u_{m+1}-\chi(\frac{|m|}{M})u_{m}]\nonumber\\
&=
  \sum\limits_{m\in
  \mathbb{Z}}\chi(\frac{|m|}{M})|(Bu)_m|^2+\sum\limits_{m\in
  \mathbb{Z}}\chi'(\frac{\tilde{m}}{M})\frac{1}{M}(u_{m+1}-u_m)u_{m+1}\nonumber\\
&\geqslant
  \sum\limits_{m\in
  \mathbb{Z}}\chi(\frac{|m|}{M})|(Bu)_m|^2-\frac{2\chi_0R^2}{M},\Hs
t\geqslant T_0,
\end{align}
where $\tilde{m}$ locates between $|m|$ and $|m+1|$. Likewise,
\begin{align}
  (Bv,Bq)
&\geqslant
  \sum\limits_{m\in
  \mathbb{Z}}\chi(\frac{|m|}{M})|(Bv)_m|^2-\frac{2\chi_0R^2}{M},\Hs t\geqslant
  T_0,\label{4.21}\\
  (BZ,BW)
&\geqslant
  \sum\limits_{m\in
  \mathbb{Z}}\chi(\frac{|m|}{M})|(BZ)_m|^2-\frac{2\chi_0R^2}{M\mu^2},\Hs t\geqslant
  T_0.\label{4.22}
\end{align}
Now
\begin{align}
  \lam\sum\limits_{m\in \mathbb{Z}}\chi(\frac{|m|}{M})a_mv_m
&\leqslant
  \frac{\lam}{2}\sum\limits_{m\in \mathbb{Z}}\chi(\frac{|m|}{M})v_m^2
 +
  \frac{\lam}{2}\sum\limits_{m\in
  \mathbb{Z}}\chi(\frac{m}{M})a_m^2,\label{4.23}\\
  k\alpha(1+\mu)\sum\limits_{m\in
  \mathbb{Z}}\chi(\frac{|m|}{M})u_mZ_m
&\leqslant
  \frac{k\alpha(1+\mu)}{2}\sum\limits_{m\in
  \mathbb{Z}}\chi(\frac{|m|}{M})(u_m^2+Z_m^2),\label{4.24}
\end{align}
and by (\textbf{H}1)-(\textbf{H}2), we have
\begin{align}
  \alpha\sum\limits_{m\in\mathbb{Z}}\chi(\frac{|m|}{M})u_mb_{1m}f_{1m}(\sigma)\!
&\leqslant
  \!\frac{\alpha(\lambda+k)}{4}\!\sum\limits_{m\in \mathbb{Z}}\chi(\frac{|m|}{M})u_m^2\!+\!\frac{\alpha
  \kappa c_1^2\pi^2}{\lam+k}\!\sum\limits_{m\in
  \mathbb{Z}}\chi(\frac{|m|}{M})b_{1m}^2,\label{4.25}\\
 \sum\limits_{m\in
  \mathbb{Z}}\chi(\frac{|m|}{M})v_mb_{2m}f_{2m}(\sigma)
&\leqslant
  \frac{\lam}{4}\sum\limits_{m\in \mathbb{Z}}\chi(\frac{|m|}{M})v_m^2
  +\frac{\kappa c_2^2\pi^2}{\lam}\sum\limits_{m\in
  \mathbb{Z}}\chi(\frac{|m|}{M})b_{2m}^2,\label{4.26}\\
\mu\alpha\sum\limits_{m\in
  \mathbb{Z}}\!\chi(\frac{|m|}{M})Z_mb_{3m}f_{3m}(\sigma)
&\leqslant\!
  \frac{\mu\alpha(\mu\lambda+k)}{4}\!\sum\limits_{m\in \mathbb{Z}}\!\chi(\frac{|m|}{M})Z_m^2
 \!+\!
  \frac{\mu\alpha \kappa c_3^2\pi^2}{\mu\lam+k}\!\sum\limits_{m\in
  \mathbb{Z}}\!\chi(\frac{|m|}{M})b_{3m}^2.\label{4.27}
\end{align}
Thus we conclude from \eqref{4.19}-\eqref{4.27} that
\begin{align*}
&
  \frac{\mathrm{d}}{\mathrm{d}t}\sum\limits_{m\in\mathbb{Z}}\chi(\frac{|m|}{M})(\alpha
  u_m^2+v_m^2+\mu^2\alpha Z_m^2)
 +
  [\frac{3}{2}\lam-(\mu-\frac{1}{2})k]\alpha\sum\limits_{m\in
  \mathbb{Z}}\chi(\frac{|m|}{M})u_m^2 \nonumber \\
 &+
  \frac{\lam}{2}\sum\limits_{m\in
  \mathbb{Z}}\chi(\frac{|m|}{M})v_m^2+[\frac{3\lam}{2}
  +
   \frac{(\frac{\mu}{2}-1)k}{\mu^2}]\mu^2\alpha\sum\limits_{m\in
  \mathbb{Z}}\chi(\frac{|m|}{M})Z_m^2 \nonumber \\
 \leqslant&\,
  \lam\sum\limits_{|m|\geqslant
  M}a_m^2+(\frac{2\alpha c_1^2}{\lam+k}\sum\limits_{|m|\geqslant
  M}b_{1m}^2+\frac{2c_2^2}{\lam}\sum_{|m|\geqslant
  M}b_{2m}^2+\frac{2\mu\alpha c_3^2}{\mu\lam+k}\sum_{|m|\geqslant
  M}b_{3m}^2)\kappa\pi^2 \nonumber \\
 &+
  \frac{4\chi_0R^2}{M}(d_1\alpha+d_2+d_3\alpha),\Hs t\geqslant T_0.
\end{align*}
Then
\begin{align}\label{4.28}
&
 \frac{\mathrm{d}}{\mathrm{d}t}\sum\limits_{m\in\mathbb{Z}}\chi(\frac{|m|}{M})(\alpha
  u_m^2+v_m^2+\mu^2\alpha Z_m^2)+\theta_1\sum\limits_{m\in\mathbb{Z}}\chi(\frac{|m|}{M})(\alpha
  u_m^2+v_m^2+\mu^2\alpha Z_m^2)\nonumber\\
\leqslant& \,
  C(M,a,b_1,b_2,b_3),\Hs t\geqslant T_0,
\end{align}
where \begin{align*}
  C(M,a,b_1,b_2,b_3)
&=
  (\frac{2\alpha c_1^2}{\lam+k}\sum\limits_{|m|\geqslant
  M}b_{1m}^2+\frac{2c_2^2}{\lam}\sum_{|m|\geqslant
  M}b_{2m}^2+\frac{2\mu\alpha c_3^2}{\mu\lam+k}\sum_{|m|\geqslant
  M}b_{3m}^2)\kappa\pi^2\\
&+
  \lam\sum\limits_{|m|\geqslant
  M}a_m^2+\frac{4\chi_0R^2}{M}(d_1\alpha+d_2+d_3\alpha).
\end{align*}
Applying Gronwall inequality to \eqref{4.28}, we obtain for any
$t\geqslant T_0,$
\begin{align*}
 &\sum\limits_{m\in\mathbb{Z}}\chi(\frac{|m|}{M})(\alpha
  u_m^2+v_m^2+\mu^2\alpha Z_m^2)=
  \sum\limits_{m\in\mathbb{Z}}\chi(\frac{|m|}{M})(\alpha
  u_m^2+v_m^2+\alpha z_m^2)\nonumber\\
\leqslant& \,
  \mathrm{e}^{-\theta_1 (t-t_0)}\sum\limits_{m\in\mathbb{Z}}\chi(\frac{|m|}{M})(\alpha
  u_{m,0}^2+v_{m,0}^2+\alpha
  z_{m,0}^2)+\frac{C(M,a,b_1,b_2,b_3)}{\theta_1}.
\end{align*}
Therefore,
\begin{align}\label{4.29}
\sum\limits_{m\in
  \mathbb{Z}}\chi(\frac{|m|}{M})|\vp_m(t)|^2
&\leqslant
  \frac{\de_2}{\de_1}\mathrm{e}^{-\theta_1(t-t_0)}\|\vp_0\|^2+\frac{C(M,a,b_1,b_2,b_3)}{\de_1\theta_1}\nonumber\\
&\leqslant
  \frac{\de_2}{\de_1}R^2\mathrm{e}^{-\theta_1(t-t_0)}+\frac{C(M,a,b_1,b_2,b_3)}{\de_1\theta_1},\Hs
  t\geqslant T_0.
\end{align}
By the assumption (\textbf{H}1), there exists $N(\varepsilon,B_0)\in
\mathbb{N}$ so that \be\label{4.30}
  \frac{C(M,d_1,d_2,d_3)}{\de_1\theta_1}
\leqslant
  \frac{\varepsilon^2}{2},\Hs \forall M\geqslant N(\varepsilon,B_0),
\ee and there is $t(\varepsilon,B_0)>T_0$ such that \be\label{4.31}
  \frac{\de_2}{\de_1}R^2\mathrm{e}^{-\theta_1 (t-t_0)}
\leqslant
  \frac{\varepsilon^2}{2}, \Hs \forall t\geqslant t(\varepsilon,B_0).\ee
Thus let $N_1(\varepsilon,B_0)=2N(\varepsilon,B_0)$ and
$t_1(\varepsilon,B_0)=t(\varepsilon,B_0)$. We deduce from
\eqref{4.29}-\eqref{4.31} that
$$
 \sum\limits_{|m|\geqslant N_1(\varepsilon,B_0)}|\vp_m(t)|^2
\leqslant
  \sum\limits_{m\in
  \mathbb{Z}}\chi(\frac{|m|}{N(\varepsilon,B_0)})|\vp_m(t)|^2
\leqslant \varepsilon^2, \Hs \forall t\geqslant t_1(\varepsilon,B_0).
$$
The proof is complete.
\eo
By Theorems \ref{t3.1}, \ref{t4.1} and Lemmas \ref{l4.2}, \ref{l4.3}, we have the
following result.
\bt\label{t4.4}
Let the assumptions $(\mathbf{H}1)$-$(\mathbf{H}2)$ hold.
Then for any $\ve>0$, there exists a family of sets $\{\cA_\ve(\sig)\}_{\sig\in\T^\kappa}$ with $$A(\sig)\subset\cA_\ve(\sig)\subset \cN_\ve(\sig),\Hs \A\sig\in \T^\kappa,$$ such that it is forward invariant for the family of processes $\{U_\sig(t,0)\}_{t\geq 0},\sig\in \T^\kappa$. Moreover, $\{\cA_\ve(\sig)\}_{\sig\in\T^\kappa}$ uniformly (w.r.t. $\sig\in \T^\kappa$) forward attracts each bounded set $B\ss E$: $$\lim_{t\ra + \infty}\sup_{\sig\in \T^\kappa} \d_H(U_\sig(t,0)B,\cA_\ve(\theta_t\sig))=0.$$
\et

\subsection{Existence of exponentially attracting set}
In this subsection, we show that the assumptions ({\bf A}1)-({\bf A}2) in Theorem \ref{t3.6} are satisfied.
\bl\label{l4.5}
Assume the assumptions $(\mathbf{H}1)$-$(\mathbf{H}2)$ hold.  Then
\benu
\item[(1)] for any $T>T_0$, there exists $L_T>0$ such that for any $\sig\in \T^\kappa$ and $t\in [T_0,T]$, $$\|U_\sig(t,0)\vp^{(1)}_0-U_\sig(t,0)\vp^{(2)}_0\|\leqslant L_T\|\vp^{(1)}_0-\vp^{(2)}_0\|,\hs \vp^{(1)}_0,\vp^{(2)}_0\in B_0;$$
\item[(2)] there exist $T^*>T_0, 0\leqslant \beta<1/2$ and a finite dimensional projection $P_{N^*}:E\mapsto E_{N^*}(N^*\in\N)$ such that for any $\sig\in \T^\kappa$, $$\|(I-P_{N^*})(U_\sig(T^*,0)\vp^{(1)}_0-U_\sig(T^*,0)\vp^{(2)}_0)\|\leqslant \beta\|\vp^{(1)}_0-\vp^{(2)}_0\|,\hs \vp^{(1)}_0,\vp^{(2)}_0\in B_0.$$  
\eenu
\el
\bo{\rm (i)} For each $\sigma \in \mathbb{T}^\kappa$, let
$$
  \varphi^{(i)}(t)=U_\sigma(t,0)\varphi^{(i)}_0=(u^{(i)}(t),v^{(i)}(t),z^{(i)}(t)),
  \hs \forall t\geqslant 0
$$
be two solutions of equations \eqref{4.5}-\eqref{4.6} with initial values
$\varphi^{(i)}_0\in B_0$ for $i=1,2$. Then, if
$t\geqslant T_0$, $\varphi^{(1)}(t),\varphi^{(2)}(t)\in B_0$ .
Write
\begin{align*}
  u_d(t)=&u^{(1)}(t)-u^{(2)}(t),\hs
  v_d(t)=v^{(1)}(t)-v^{(2)}(t),\nonumber\\
  z_d(t)=&z^{(1)}(t)-z^{(2)}(t),\hs
  \varphi_d(t)=\varphi^{(1)}(t)-\varphi^{(2)}(t).
\end{align*}
By equations \eqref{4.5}-\eqref{4.6}, we have
\begin{align}
&\dot{\varphi}_d+\Theta\varphi_d=G(\varphi^{(1)},\sigma,t)-G(\varphi^{(2)},\sigma,t),\label{4.32}\\
 & \varphi_d(0)=\varphi^{(1)}_0-\varphi^{(2)}_0.\nonumber
\end{align}
Taking the inner product of \eqref{4.32} with $\vp_d$ in $E$, one obtains
\be\label{4.33}
\frac{1}{2}\frac{\mathrm{d}}{\mathrm{d}t}\|\vp_{d}\|^2+\big(\Theta\varphi_d-G(\varphi^{(1)},\sigma,t)+G(\varphi^{(2)},\sigma,t),\vp_d\big)=0.
\ee
Since $\Theta: E\mapsto E$ is a bounded linear operator, $G:\~E\X\R\mapsto E$ is locally Lipschitz continuous with respect to $\vp$ for each $\sig\in\T^\kappa$ and $t\in\R$ and the set $B_0$ is bounded, we conclude that there exist positive constants $C_0$ and $L_0$ such that
\begin{align}
&\hs \big(\Theta\varphi_d-G(\varphi^{(1)},\sigma,t)+G(\varphi^{(2)},\sigma,t),\vp_d\big)\nonumber\\
&\leqslant
\big(C_0\|\vp_d\|+\|G(\varphi^{(2)},\sig,t)-G(\varphi^{(1)},\sig,t)\|\big)\|\vp_d\|\nonumber\\
&\leqslant
(C_0+L_0)\|\vp_d\|^2.\label{4.34}
\end{align}
Thus it follows from \eqref{4.33}-\eqref{4.34} that
\be\label{4.35}
\frac{\mathrm{d}}{\mathrm{d}t}\|\vp_{d}\|^2\leqslant C_1\|\vp_d\|^2,
\ee
where $C_1=2(C_0+L_0)$. Applying Gronwall inequality to \eqref{4.35} on $[0,t]$ with $t\in [T_0,T]$, we have
\be\label{4.36}
\|\vp_d(t)\|^2\leqslant {\rm e}^{C_1t}\|\vp_d(0)\|^2,
\ee
which implies that
\begin{align*}
  \|\vp^{(1)}(t)-\vp^{(2)}(t)\|&=\|U_\sig(t,0)\vp_0^{(1)}-U_\sig(t,0)\vp_0^{(2)}\|\nonumber \\
&\leqslant
  L_T\|\vp_0^{(1)}-\vp_0^{(2)}\|,\Hs t\in [T_0,T],
\end{align*}
where $L_T=\sqrt {{\rm e}^{C_1T}}.$
\vs
{\rm (ii)} Set
\begin{align*}
  p_d&=(p_{dm})_{m\in \mathbb{Z}},\ q_d=(q_{dm})_{m\in \mathbb{Z}},\  w_d=(w_{dm})_{m\in \mathbb{Z}}, \\
  p_{dm}&=\chi(\frac{|m|}{M})u_{dm},\ q_{dm}=\chi(\frac{|m|}{M})v_{dm},\, w_{dm}=\chi(\frac{|m|}{M})z_{dm},
\end{align*}
where $M$ is a positive integer which will be decided later and
$\chi(x)$ is defined by \eqref{4.18}. By equation \eqref{4.1}, we have
\begin{align}\label{4.37}
 \dot{u}_d
=&
  -d_1 Au_{d}-(\lambda+k)u_d+(u^{(1)})^2v^{(1)}-(u^{(2)})^2v^{(2)}-\alpha
  \((u^{(1)})^3-(u^{(2)})^3\) \nonumber \\
&+\beta z_d.
\end{align}
Taking the inner product of \eqref{4.37} with $p_d$ in $\ell^2$ gives
\begin{align}\label{4.38}
  \frac{1}{2}\frac{\mathrm{d}}{\mathrm{d}t}\sum\limits_{m\in{\bf\mathbb{Z}}}\chi(\frac{|m|}{M})u_{dm}^2
&=
  \!-d_1(Au_{d},p_d)\!-\!(\lambda+k)(u_d,p_d)\!+\!\((u^{(1)})^2v^{(1)}-(u^{(2)})^2v^{(2)},p_d\)\nonumber\\
&
  -\alpha\((u^{(1)})^3-(u^{(2)})^3,p_d\)+\beta (z_d,p_d).
\end{align}
By some computations, we have
\begin{align}\label{4.39}
  (Au_{d},p_d)
&=
  (Bu_{d},Bp_d)
 =
  \sum\limits_{m\in{\bf\mathbb{Z}}}(Bu_d)_m(\chi(\frac{|m+1|}{M})u_{dm+1}-\chi(\frac{|m|}{M})u_{dm})\nonumber\\
&=
  \sum\limits_{m\in{\bf\mathbb{Z}}}\chi(\frac{|m|}{M})(Bu_d)_m^2
 +
  \sum\limits_{m\in{\bf\mathbb{Z}}}\chi'(\frac{\tilde{m}}{M})\frac{1}{M}(u_{dm+1}-u_{dm})u_{dm+1}\nonumber\\
&\geqslant
  \sum\limits_{m\in{\bf\mathbb{Z}}}\chi(\frac{|m|}{M})(Bu_d)_m^2-\frac{2\chi_0}{M}\|\varphi_d\|^2_E,
\end{align}
where $\tilde{m}$ is a constant locating between $|m|$ and $|m+1|$.
Thanks to Lemma \ref{l4.3}, there exist $t_1 :=t_1(\alpha,\lambda, B_0)>T_0$
and $N_1 :=N_1(\alpha,\lambda,B_0)\in \mathbb{N}$ such that if
$t\geqslant t_1$ and $M\geqslant N_1$,
\begin{align}\label{4.40}
&
  \((u^{(1)})^2v^{(1)}-(u^{(2)})^2v^{(2)},p_d\)\nonumber\\
=&
  \sum\limits_{m\in{\bf\mathbb{Z}}}\chi(\frac{|m|}{M})((u^{(1)}_m)^2v^{(1)}_m-(u^{(2)}_m)^2v^{(2)}_m)u_{dm}\nonumber\\
=&
  \sum\limits_{m\in{\bf\mathbb{Z}}}\chi(\frac{|m|}{M})(u^{(1)}_m)^2v_{dm}u_{dm}
 +
  \sum\limits_{m\in{\bf\mathbb{Z}}}\chi(\frac{|m|}{M})(u_m^{(1)}+u_m^{(2)})v_m^{(2)}u_{dm}^2\nonumber\\
\leqslant&
  \frac{\lambda}{10}\sum\limits_{m\in{\bf\mathbb{Z}}}\chi(\frac{|m|}{M})|v_{dm}u_{dm}|
 +
  \frac{\lambda}{10}\sum\limits_{m\in{\bf\mathbb{Z}}}\chi(\frac{|m|}{M})u_{dm}^2\nonumber\\
\leqslant&
  \frac{\lambda}{5}\sum\limits_{m\in{\bf\mathbb{Z}}}\chi(\frac{|m|}{M})u_{dm}^2
 +
  \frac{\lambda}{40}\sum\limits_{m\in{\bf\mathbb{Z}}}\chi(\frac{|m|}{M})v_{dm}^2
\end{align}
and
\begin{align}\label{4.41}
  \((u^{(1)})^3-(u^{(2)})^3,p_d\)
&=
  \sum\limits_{m\in{\bf\mathbb{Z}}}\chi(\frac{|m|}{M})((u^{(1)}_m)^3-(u^{(2)}_m)^3)u_{dm}\nonumber\\
&\leqslant
  \sum\limits_{m\in{\bf\mathbb{Z}}}\chi(\frac{|m|}{M})(|u_m^{(1)}|+|u_m^{(2)}|)^2u_{dm}^2\nonumber\\
&\leqslant
  \frac{\lambda}{10\alpha}\sum\limits_{m\in{\bf\mathbb{Z}}}\chi(\frac{|m|}{M})u_{dm}^2.
\end{align}
Note that
\begin{align}\label{4.42}
  \beta(z_d,p_d)
 =
  \beta\sum\limits_{m\in{\bf\mathbb{Z}}}\chi(\frac{|m|}{M})z_{dm}u_{dm}
\leqslant
  \frac{\beta}{2}\sum\limits_{m\in{\bf\mathbb{Z}}}\chi(\frac{|m|}{M})z_{dm}^2
 +
  \frac{\beta}{2}\sum\limits_{m\in{\bf\mathbb{Z}}}\chi(\frac{|m|}{M})u_{dm}^2.
\end{align}
Thus taking \eqref{4.38}-\eqref{4.42} into account, we obtain
\begin{align}\label{4.43}
&\hs
  \frac{1}{2}\frac{\mathrm{d}}{\mathrm{d}t}\sum\limits_{m\in{\bf\mathbb{Z}}}\chi(\frac{|m|}{M})u_{dm}^2
 +
  (\lambda+k)\sum\limits_{m\in{\bf\mathbb{Z}}}\chi(\frac{|m|}{M})u_{dm}^2\nonumber\\
&\leqslant
  (\frac{3\lambda}{10}+\frac{\beta}{2})\sum\limits_{m\in{\bf\mathbb{Z}}}\chi(\frac{|m|}{M})u_{dm}^2
 +
  \frac{\lambda}{40}\sum\limits_{m\in{\bf\mathbb{Z}}}\chi(\frac{|m|}{M})v_{dm}^2
+
  \frac{\beta}{2}\sum\limits_{m\in{\bf\mathbb{Z}}}\chi(\frac{|m|}{M})z_{dm}^2\nonumber\\
&\quad
 +
  \frac{2d_1\chi_0}{M}\|\varphi_d\|^2_E,\Hs t\geqslant t_1,\,\,M\geqslant N_1.
\end{align}
From equation \eqref{4.2}, we have
\begin{align}\label{4.44}
 \dot{v}_d
=
  -d_2 Av_{d}-\lambda v_d-\((u^{(1)})^2v^{(1)}-(u^{(2)})^2v^{(2)}\)+\alpha
  \((u^{(1)})^3-(u^{(2)})^3\).
\end{align}
Taking the inner product of \eqref{4.44} with $q_d$ in $\ell^2$ gives
\begin{align}\label{4.45}
  \frac{1}{2}\frac{\mathrm{d}}{\mathrm{d}t}\sum\limits_{m\in{\bf\mathbb{Z}}}\chi(\frac{|m|}{M})v_{dm}^2
&=
  -d_2(Av_{d},q_d)-\lambda(v_d,q_d)-\((u^{(1)})^2v^{(1)}-(u^{(2)})^2v^{(2)},q_d\)\nonumber\\
& \hs
  +\alpha\((u^{(1)})^3-(u^{(2)})^3,q_d\).
\end{align}
Similar to \eqref{4.39}, we have
\begin{align}\label{4.46}
  (Av_{d},q_d)
\geqslant
  \sum\limits_{m\in{\bf\mathbb{Z}}}\chi(\frac{|m|}{M})(Bv_d)_m^2-\frac{2\chi_0}{M}\|\varphi_d\|^2_E.
\end{align}
By virtue of Lemma \ref{l4.3}, there exist $t_2
:=t_2(\alpha,\lambda, B_0)$ and $N_2 :=N_2(\alpha,\lambda,B_0)\in
\mathbb{N}$ with $t_2>t_1$ and $N_2>N_1$ such that if $t\geqslant t_2$
and $M\geqslant N_2$,
\begin{align}\label{4.47}
&
  \((u^{(1)})^2v^{(1)}-(u^{(2)})^2v^{(2)},q_d\) \nonumber\\
=&
  \sum\limits_{m\in{\bf\mathbb{Z}}}\chi(\frac{|m|}{M})((u^{(1)}_m)^2v^{(1)}_m-(u^{(2)}_m)^2v^{(2)}_m)v_{dm}\nonumber\\
=&
  \sum\limits_{m\in{\bf\mathbb{Z}}}\chi(\frac{|m|}{M})(u^{(1)}_m)^2v_{dm}^2
 +
  \sum\limits_{m\in{\bf\mathbb{Z}}}\chi(\frac{|m|}{M})(u_m^{(1)}+u_m^{(2)})v_m^{(2)}u_{dm}v_{dm}\nonumber\\
\leqslant&
  \frac{\lambda}{10}\sum\limits_{m\in{\bf\mathbb{Z}}}\chi(\frac{|m|}{M})u_{dm}^2
 +
  \frac{\lambda}{20}\sum\limits_{m\in{\bf\mathbb{Z}}}\chi(\frac{|m|}{M})v_{dm}^2
\end{align}
and
\begin{align}\label{4.48}
  \((u^{(1)})^3-(u^{(2)})^3,q_d\)
&=
  \sum\limits_{m\in{\bf\mathbb{Z}}}\chi(\frac{|m|}{M})((u^{(1)}_m)^3-(u^{(2)}_m)^3)v_{dm}\nonumber\\
&\leqslant
  \sum\limits_{m\in{\bf\mathbb{Z}}}\chi(\frac{|m|}{M})(|u_m^{(1)}|+|u_m^{(2)}|)^2|u_{dm}v_{dm}|\nonumber\\
&\leqslant
  \frac{\lambda}{10\alpha}\sum\limits_{m\in{\bf\mathbb{Z}}}\chi(\frac{|m|}{M})u_{dm}^2
  +
  \frac{\lambda}{40\alpha}\sum\limits_{m\in{\bf\mathbb{Z}}}\chi(\frac{|m|}{M})v_{dm}^2.
\end{align}
Combining \eqref{4.45}-\eqref{4.48}, we see that for any $t\geqslant t_2$
and $M\geqslant N_2$,
\begin{align}\label{4.49}
& \quad \,
  \frac{1}{2}\frac{\mathrm{d}}{\mathrm{d}t}\sum\limits_{m\in{\bf\mathbb{Z}}}\chi(\frac{|m|}{M})v_{dm}^2
 +
  \frac{37\lambda}{40}\sum\limits_{m\in{\bf\mathbb{Z}}}\chi(\frac{|m|}{M})v_{dm}^2 \nonumber\\
&\leqslant
 \frac{\lambda}{5}\sum\limits_{m\in{\bf\mathbb{Z}}}\chi(\frac{|m|}{M})u_{dm}^2
 +
  \frac{2d_2\chi_0}{M}\|\varphi_d\|^2_E.
\end{align}
Finally, by \eqref{4.3}, we have
\begin{align}\label{4.50}
 \dot{z}_d
=
  -d_3 Az_{d}+k u_d-(\lambda+\beta)z_d.
\end{align}
Taking the inner product of \eqref{4.50} with $w_d$ in $\ell^2$ gives
\begin{align}\label{4.51}
  \frac{1}{2}\frac{\mathrm{d}}{\mathrm{d}t}\sum\limits_{m\in{\bf\mathbb{Z}}}\chi(\frac{|m|}{M})z_{dm}^2
=&
  -d_3(Az_{d},w_d)+k(u_d,w_d)-(\lambda+\beta)(z_d,w_d).
\end{align}
Similar to the derivation of \eqref{4.49}, we deduce that
\begin{align}\label{4.52}
& \quad \,
  \frac{1}{2}\frac{\mathrm{d}}{\mathrm{d}t}\sum\limits_{m\in{\bf\mathbb{Z}}}\chi(\frac{|m|}{M})z_{dm}^2
 +
  (\lam+\beta-\frac{k}{2})\sum\limits_{m\in{\bf\mathbb{Z}}}\chi(\frac{|m|}{M})z_{dm}^2\nonumber\\
&\leqslant
  \frac{k}{2}\sum\limits_{m\in{\bf\mathbb{Z}}}\chi(\frac{|m|}{M})u_{dm}^2
 +
  \frac{2d_3\chi_0}{M}\|\varphi_d\|^2_E.
\end{align}
Thus, we conclude from \eqref{4.43}, \eqref{4.49} and \eqref{4.52}
that if $t\geqslant t_2$ and $M\geqslant N_2$,
\begin{align*}
&
  \frac{1}{2}\frac{\mathrm{d}}{\mathrm{d}t}\sum\limits_{m\in{\bf\mathbb{Z}}}\chi(\frac{|m|}{M})(u_{dm}^2+v_{dm}^2+z_{dm}^2)
 +
 \frac{\lam+k-\beta}{2}\sum\limits_{m\in{\bf\mathbb{Z}}}\chi(\frac{|m|}{M})u_{dm}^2
 +
  \frac{9\lambda}{10}\sum\limits_{m\in{\bf\mathbb{Z}}}\chi(\frac{|m|}{M})v_{dm}^2\nonumber\\
&+
  \frac{\lambda+\beta-k}{2}\sum\limits_{m\in{\bf\mathbb{Z}}}\chi(\frac{|m|}{M})z_{dm}^2
\leqslant
  \frac{2\chi_0(d_1+d_2+d_3)}{M}\|\varphi_d\|^2_E,
\end{align*}
which implies that for any $t\geqslant t_2$ and $M\geqslant N_2$,
\begin{align}\label{4.53}
&
  \frac{\mathrm{d}}{\mathrm{d}t}\sum\limits_{m\in{\bf\mathbb{Z}}}\chi(\frac{|m|}{M})|\varphi_{dm}|_{E}^2
 +
  \theta_2\sum\limits_{m\in{\bf\mathbb{Z}}}\chi(\frac{|m|}{M})|\varphi_{dm}|_{E}^2
\leqslant
  \frac{4\chi_0(d_1+d_2+d_3)}{M}\|\varphi_d\|^2_E,
\end{align}
where
$\theta_2=\min\{\lam+k-\beta,\frac{9\lambda}{5},\lambda+\beta-k\}$.
Applying Gronwall inequality to \eqref{4.53} from $t_2$ to $t$ with
$t\geqslant t_2$, we have
\begin{align}\label{4.54}
  \sum\limits_{m\in{\bf\mathbb{Z}}}\chi(\frac{|m|}{M})|\varphi_{dm}(t)|_{E}^2
\leqslant&\,
  \mathrm{e}^{-\theta_2(t-t_2)}\sum\limits_{m\in{\bf\mathbb{Z}}}\chi(\frac{|m|}{M})|\varphi_{dm}(t_2)|_{E}^2\nonumber\\
&+
  \frac{4\chi_0(d_1+d_2+d_3)}{M}\int_{t_2}^{t}\|\varphi_d(s)\|_E^2\mathrm{e}^{-\theta_2(t-s)}\mathrm{d}s,
\end{align}
provided $M\geqslant N_2.$
By \eqref{4.36}, one has
\begin{align}\label{4.55}
\mathrm{e}^{-\theta_2(t-t_2)}\sum\limits_{m\in{\bf\mathbb{Z}}}\chi(\frac{|m|}{M})|\varphi_{dm}(t_2)|_{E}^2\leqslant \mathrm{e}^{-\theta_2(t-t_2)}\|\vp_d(t_2)\|_E^2\leqslant \mathrm{e}^{-\theta_2(t-t_2)+C_1t_2}\|\vp_d(0)\|_E^2
\end{align}
and
\begin{align}\label{4.56}
&\frac{4\chi_0(d_1+d_2+d_3)}{M}\int_{t_2}^{t}\|\varphi_d(s)\|_E^2\mathrm{e}^{-\theta_2(t-s)}\mathrm{d}s\nonumber\\
\leqslant&
\frac{4\chi_0(d_1+d_2+d_3)}{M}\|\varphi_d(0)\|_E^2\int_{t_2}^{t}\mathrm{e}^{-\theta_2(t-s)+C_1s}\mathrm{d}s\nonumber\\
\leqslant&
\frac{4\chi_0(d_1+d_2+d_3)}{M(\theta_2+C_1)}\|\varphi_d(0)\|_E^2\mathrm{e}^{C_1t}.
\end{align}
Thus it follows from \eqref{4.54}-\eqref{4.56} that if $t\geqslant t_2$ and $M\geqslant N_2$,
\begin{align}\label{4.57}
\hs \sum\limits_{m\in{\bf\mathbb{Z}}}\chi(\frac{|m|}{M})|\varphi_{dm}(t)|_{E}^2
\leqslant& \,\mathrm{e}^{-\theta_2(t-t_2)+C_1t_2}\|\vp_d(0)\|_E^2\nonumber\\
&\,+\frac{4\chi_0(d_1+d_2+d_3)}{M(\theta_2+C_1)}\|\varphi_d(0)\|_E^2\mathrm{e}^{C_1t}.
\end{align}
Now pick two constants $T^*\geqslant t_2$ and $N_3\geqslant N_2$ such that
$$
 \mathrm{e}^{-\theta_2(T^*-t_2)+C_1t_2}
 +\frac{4\chi_0(d_1+d_2+d_3)}{N_3(\theta_2+C_1)}\mathrm{e}^{C_1T^*}:=\beta^2<\frac{1}{4}.$$
Therefore if $N^*\geqslant 2N_3$, we deduce from \eqref{4.57} that
$$
\sum\limits_{|m|\geqslant N^*}|\varphi_{dm}(T^*)|_{E}^2\leqslant \sum\limits_{m\in{\bf\mathbb{Z}}}\chi(\frac{|m|}{N_3})|\varphi_{dm}(T^*)|_{E}^2
\leqslant \beta^2\|\varphi_d(0)\|_E^2,
$$
which implies
$$\|(I-P_{N^*})(U_\sig(T^*,0)\vp_0^{(1)}-U_\sig(T^*,0)\vp_0^{(2)})\|\leqslant \beta\|\vp_0^{(1)}-\vp_0^{(2)}\|,$$  %
where $\beta< 1/2.$ This completes the proof of the lemma.
\eo

According to Lemmas \ref{l4.2}, \ref{l4.5} and Theorem \ref{t3.6}, we have
\bt\label{t4.6}
Assume the assumptions $(\mathbf{H}1)$-$(\mathbf{H}2)$ hold.  Then for each $\ve>0$, there exists a family of sets $\{\cB_\ve(\sig)\}_{\sig\in \T^\kappa}$ which is forward invariant for the family of processes $\{U_\sig(t,0)\}_{t\geq 0},\sig\in \T^\kappa$ generated by equations \eqref{4.5}-\eqref{4.6} such that  $\{\cB_\ve(\sig)\}_{\sig\in \T^\kappa}$ uniformly forward exponentially attracts each bounded subset of $E$.
\et
\setcounter {equation}{0}
\section{Remarks on some extensions}
We end the article with two remarks on some possible extensions of our main results.
\br
Our main results could be examined for the general lattice system \eqref{e3.1}-\eqref{e3.2} considered in $\Z^n$ for some positive integer $n\geqslant 2$,
\begin{align}
&\dot {u}=f(u)+\sig(t),\hs u=(u_m)_{m\in\Z^n}, \,\,t>t_0;\label{5.1}\\
&u(t_0)=u_0,\hs u_0=(u_{m,0})_{m\in\Z}\in l^2, \label{5.2}
\end{align}
where $$l^2=\{u=(u_{m})_{m\in \mathbb{Z}^n}:\  u_{m}\in\mathbb{R},\hs
  \sum\limits_{m\in\mathbb{Z}^n}u_{m}^{2}<+\infty\}.
$$
Moreover, the corresponding results of Theorem \ref{t4.4} and Theorem \ref{t4.6} could be verified for the following LDS
with quasi-periodic external forces:
 \begin{align}
\dot{u}_{m}
&=-d_1(Au)_m-(\lambda+k)u_m+u_m^2v_m-\alpha
u_m^3+\beta z_m+b_{1m}f_{1m}(\sigma(t)),\label{5.3}\\
  \dot{v}_{m}
&=-d_2(Av)_m+\lambda(a_m-v_m)-u_m^2v_m+\alpha
u_m^3+b_{2m}f_{2m}(\sigma(t)), \label{5.4}\\
  \dot{z}_{m}
&=-d_3(Az)_m+ku_m-(\lambda+\beta)z_m+b_{3m}f_{3m}(\sigma(t))\label{5.5}
 \end{align}
for $m=(m_1,m_2,\cdot\.\.,m_n)\in \mathbb{Z}^n,$ where the operator $A$ is defined as
\begin{align*}
(Au)_{m}&=(Au)_{(m_1,m_2,\.\.\.,m_n)}\\
&=2ku_m-u_{(m_1+1,m_{2},\.\.\.,m_n)}-u_{(m_1,m_2+1,\.\.\.,m_n)}-\.\.\.-u_{(m_1,m_2,\.\.\.,m_n+1)}\\
&\hs -u_{(m_1-1,m_{2},\.\.\.,m_n)}-u_{(m_1,m_2-1,\.\.\.,m_n)}-\.\.\.-u_{(m_1,m_2,\.\.\.,m_n-1)}.
\end{align*}
Thus the similar results of Theorem \ref{t4.4} and Theorem \ref{t4.6} for equations \eqref{5.3}-\eqref{5.5} are still valid. In such case, equations \eqref{5.3}-\eqref{5.5} can be regarded as a discrete analogue of nonautonomous Gray-Scott equations \eqref{aa6}-\eqref{aa8} in $\R^n$.
\er
\br
In this article, we consider the forward dynamical behavior of nonautonomous LDSs of the form
$\dot {u}=f(u)+\sig(t)$. In fact, for a nonautonomous LDS:
$\dot {u}=f(u,t),$
if the symbol space is compact and the corresponding conditions are satisfied, one can also obtain the similar results. We will investigate this issue in another paper.
\er
\Vs

\medskip
\medskip

\end{document}